\documentclass[11pt]{article}

\usepackage{amsmath,amssymb,latexsym}
\usepackage{theorem}
\usepackage[final]{graphicx}
\newtheorem{theorem}{Theorem}
\newtheorem{lemma}{Lemma}

\newtheorem{proposition}{Proposition}
{\theorembodyfont{\rmfamily} 

\newtheorem{remark}{Remark}}
{\theorembodyfont{\slshape} }

\newcommand{\res}{\mathop{\rm res}}

\newcommand{\supp}{\mathop{\rm supp}}
\newcommand{\field}[1]{\mathbb{#1}}
\newcommand{\R}{\field{R}}

\newcommand{\N}{\field{N}}
\newcommand{\C}{\field{C}}

\newcommand{\CC}{{\mathcal C}}

\newcommand{\NN}{{\mathcal N}}

\newcommand{\Pc}{\mathop{\rm Pc}}
\newcommand{\Int}{\mathop{\rm Int}}
\newcommand{\Ai}{\operatorname{\rm  Ai}}

\renewcommand{\Re}{\mathop{\rm Re}}
\renewcommand{\Im}{\mathop{\rm Im}}
\newcommand{\dist}{\mathop{\rm dist}}
\newcommand{\isdef}{\stackrel{\text{\tiny def}}{=}}

{\rm \trivlist \item[\hskip \labelsep{\bf Proof. }]}%
{\hspace*{\fill}$\Box$\endtrivlist}
\newenvironment{varproof}%
{\rm \trivlist \item[\hskip \labelsep{\bf Proof}]}%
{\hspace*{\fill}$\Box$\endtrivlist}
\numberwithin{equation}{section}

\title{Riemann-Hilbert analysis for Jacobi polynomials orthogonal
on a single contour}

\author{A.\ Mart\'{\i}nez-Finkelshtein\\ University of Almer\'{\i}a and Instituto
Carlos I de F\'{\i}sica Te\'{o}rica \\ y Computacional, Granada
University (SPAIN), \\ E-mail: {\tt andrei@ual.es}  \and R.\
Orive\footnote{Corresponding author.}
\\ University of La Laguna, Canary Islands (SPAIN), \\ E-mail: {\tt
rorive@ull.es}
 }

\begin{document}

\pagestyle{myheadings} \thispagestyle{plain}
\markboth{MART\'{I}NEZ-FINKELSHTEIN AND ORIVE}{RIEMANN-HILBERT
ANALYSIS FOR JACOBI POLYNOMIALS} \maketitle

\begin{abstract} Classical Jacobi polynomials
$P_{n}^{(\alpha,\beta)}$, with $\alpha, \beta>-1$, have a number
of well-known properties, in particular the location of their
zeros in the open interval $(-1,1)$. This property is no longer
valid for other values of the parameters; in general, zeros are
complex. In this paper we study the strong asymptotics of Jacobi
polynomials where the real parameters $\alpha_n,\beta_n$ depend on
$n$ in such a way that
$$
\lim_{n\rightarrow\infty}\frac{\alpha_{n}}{n}=A\,, \quad
\lim_{n\rightarrow\infty}\frac{\beta_{n}}{n}=B\,,
$$
with $A,B \in \mathbb{R}$. We restrict our attention to the case
where the limits $A,B$ are not both positive and take values
outside of the triangle bounded by the straight lines $A=0$, $B=0$
and $A+B+2=0$. As a corollary, we show that in the limit the zeros
distribute along certain curves that constitute trajectories of a
quadratic differential.

The non-hermitian orthogonality relations for Jacobi polynomials
with varying parameters lie in the core of our approach; in the
cases we consider, these relations hold on a single contour of the
complex plane. The asymptotic analysis is performed using the
Deift-Zhou steepest descent method based on the Riemann-Hilbert
reformulation of Jacobi polynomials.
\end{abstract}

\section{Introduction}
Jacobi polynomials $P_n^{(\alpha,\beta)}$ are given explicitly by
\begin{equation}\label{explJac}
P_n^{(\alpha,\beta)} (z)=2^{-n} \sum_{k=0}^n \binom{n+\alpha}{n-k}
\binom{n+\beta}{k}(z-1)^k (z+1)^{n-k}\,,
\end{equation}
or, equivalently, by the Rodrigues formula \cite[Ch.\
IV]{szego:1975}
\begin{equation}\label{RodrJac}
P_n^{(\alpha,\beta)} (z)=\frac{1}{2^n n!} (z-1)^{-\alpha}
(z+1)^{-\beta} \left( \frac{d}{dz} \right)^n \left[
(z-1)^{n+\alpha} (z+1)^{n+\beta}\right]\,.
\end{equation}
In the classical situation ($\alpha, \beta>-1$) the  Jacobi
polynomials are orthogonal in $[-1,1]$ with respect to the weight
function $(1-x)^\alpha (1+x)^\beta$ and, consequently, their zeros
are simple and located in $(-1,1)$.

Expressions (\ref{explJac})--(\ref{RodrJac}) show that the
definition of $P_n^{(\alpha,\beta)}$ may be extended to arbitrary
$\alpha,\beta \in \mathbb{R}$ (or even $\C$); but some properties
of the classical polynomials, in particular the location and
simplicity of the zeros, are no longer valid. In fact,
$P_{n}^{\left( \alpha ,\beta \right) }$ may have a multiple zero
at $z=1$ if $\alpha \in \{-1,\ldots,-n\}$, at $z=-1$ if $\beta \in
\{ -1,\ldots,-n\} $ or at $z=\infty $ (which means a degree
reduction) if $n+\alpha +\beta \in \{-1, \ldots,-n\}$.

More precisely, for $k\in \{1,\ldots,n\}$, we have (see
\cite[formula (4.22.2)]{szego:1975}),
\begin{equation}
P_{n}^{\left( -k,\beta \right) }(z) = \frac{\Gamma(n+\beta+1)
}{\Gamma(n+\beta +1-k)}\, \frac{(n-k)!}{n!} \left(
\frac{z-1}{2}\right) ^{k}P_{n-k}^{\left( k,\beta \right) }\left(
z\right). \label{integer 1}
\end{equation}
This implies in particular that $P_{n}^{( -k,\beta) }(z) \equiv 0$
if additionally $\max \left\{ k,-\beta \right\} \leq n\leq k-\beta
-1$. Analogous relations hold for $P_{n}^{\left( \alpha ,-l\right)
}$ when $l\in \{1,\ldots,n\}$. Thus, when both $k,l\in \mathbb{N}$
and $k+l\leq n $, we have
\begin{equation}
P_{n}^{\left( -k,-l\right) }\left( z\right) =2^{-k-l}\left(
z-1\right) ^{k}\left( z+1\right) ^{l}P_{n-k-l}^{\left( k,l\right)
}\left( z\right)\,. \label{integer 2}
\end{equation}
Furthermore, when $n+\alpha +\beta  =-k \in \{-1, \ldots,-n\}$,
\begin{equation}
P_{n}^{\left( \alpha ,\beta  \right) }\left( z\right)
=\frac{\Gamma(n+\alpha+1) }{\Gamma(k+\alpha )} \,
\frac{(k-1)!}{n!}\, P_{k-1}^{\left( \alpha ,\beta \right) }\left(
z\right), \label{integer 3}
\end{equation}
see \cite[Eq.\ (4.22.3)]{szego:1975}; we refer the reader to
\cite[\S 4.22]{szego:1975} for a more detailed discussion. Taking
into account formulas (\ref{integer 1})--(\ref{integer 3}) we
exclude these special integer parameters from our further
analysis.

\begin{figure}[htb] \label{fig:Jacobi_cases}
\centering \includegraphics[scale=0.9]{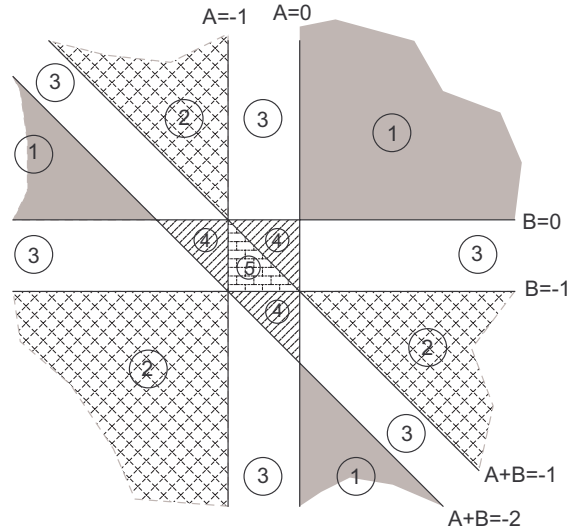} \caption{The
five different cases in the classification of Jacobi polynomials
with varying parameters according to \cite{MR2002d:33017}.}
\end{figure}

In this paper we study the asymptotic behavior of Jacobi
polynomials $P_n^{(\alpha_n,\beta_n)}$, where the parameters
$\alpha_n,\beta_n$ depend on the degree $n$ in such a way that
\begin{equation}\label{AB}
    \lim_{n\rightarrow \infty}\frac{\alpha_n}{n}\,=\,A\,, \quad \lim_{n\rightarrow
    \infty}\frac{\beta_n}{n}=B\,,
\end{equation}
and
\begin{equation}\label{region}
     A<0<B,\quad A+B>-1, \quad A \neq -1.
\end{equation}
In \cite{MR2002d:33017}, the authors considered different regions
of the $(A,B)$-plane, corresponding to different cases in the
asymptotic study of Jacobi polynomials with varying parameters
satisfying (\ref{AB}). The symmetry relations (see \cite[Ch.
IV]{szego:1975})
\begin{equation}\label{symmetry 1}
P_n^{(\alpha,\beta)}(z)\,=\,(-1)^{n}P_n^{(\beta,\alpha)}(z)
\end{equation}
and
\begin{equation}\label{symmetry 2}
P_n^{(\alpha,\beta)}(z)\,= \,\left(\frac{1-z}{2}\right)^n
P_n^{(\alpha',\beta)}\left(\frac{z+3}{z-1}\right),
\end{equation}
where $\alpha'=-2n-\alpha- \beta -1$, allow to restrict our study
to the following cases, from which all the others can be obtained:
\begin{gather}
\tag{C.1} A,B>0\,; \label{c1}\\ \tag{C.2} A<-1,\,\quad \text{and} \quad A + B >-1\,; \label{c2}\\
\tag{C.3}  -1 < A< 0 \quad \text{and} \quad B>0 \,; \label{c3}
\\ \tag{C.4}   A+B>-1,\, \quad  \text{and} \quad  A,B<0 \,; \label{c4}
\\ \tag{C.5}   A+B<-1 \quad \text{and} \quad  A,B>-1 \, \label{c5}
\end{gather}
(see Figure \ref{fig:Jacobi_cases}, which appeared first in
\cite{MR2002d:33017}, where equivalent regions under those
transformations are shown).

Case C.1 is classical and has been widely studied (see
\cite{Bosbach99}, \cite{ChenIsmail}, \cite{DetteStudden:95},
\cite{Gawronski91}, \cite{KuijlaarsVanAssche:99} and
\cite{MoakSaffVarga}). The asymptotic results therein are based on
either the well known orthogonality conditions satisfied by the
Jacobi polynomials on $[-1,1]$, or on their integral
representation.

However, until very recently, strong asymptotics of
$P_n^{(\alpha_n, \beta_n ) }$, when $\alpha _n, \beta _n$ take
\emph{arbitrary real} values and limits \eqref{AB} exist, has not
been established. In this case the orthogonality conditions were
unknown and the complex saddle points make the application of the
classical steepest descent method to the integral representation
of Jacobi polynomials practically unfeasible.

A non-hermitian orthogonality satisfied by Jacobi polynomials in
the case C.2 has been observed in \cite{MR2002d:33017}; this fact
was used there to establish the asymptotic zero distribution using
a potential theory approach.

Recently in \cite{Etna} a whole spectrum of orthogonality
conditions for Jacobi polynomials with arbitrary real parameters
has been established. In particular, we can find examples of
orthogonality on a contour or arc of the complex plane, an
incomplete or quasi-orthogonality, or even multiple or
Hermite-Pad\'{e} orthogonality conditions. The classification of the
cases depends on the number of inequalities $-1 < A < 0$, $-1 < B
< 0$, $-2 < A + B < -1$ that are satisfied. In particular, cases
C.3, C.4 and C.5 correspond to combination of parameters $A$ and
$B$ such that exactly one, exactly two, or exactly three,
respectively, of the inequalities are satisfied (cf. Figure
\ref{fig:Jacobi_cases}).

Nevertheless, the method used in \cite{MR2002d:33017} cannot be
immediately extended to the rest of the cases. One of the
essential assumptions there is a non-hermitian orthogonality of
the polynomials on a single contour, on which the support of a
certain equilibrium measure has a connected complement.

Due to this reason, in \cite{ArnoAndrei} the steepest descent
method of Deift and Zhou (see \cite{DeiftZhou}), based on a matrix
Riemann-Hilbert problem, was used to establish the strong uniform
asymptotics of the Jacobi polynomials with parameters satisfying
conditions C.5. In this paper, we use several results and ideas
from there.

The aim of the present article is to extend this analysis to
sequences of Jacobi polynomials with varying parameters
corresponding to cases C.2--C.3. Note that along with the case
C.1, these are the only situations when a full system of
orthogonality relations on a single contour in $\C$ exists.

We also remark  that a similar study, but for Laguerre polynomials
with varying parameters, has been carried out in
\cite{Kuijlaars01}-\cite{Arnoken2} and \cite{MR1858305}.

The paper is organized as follows. In Section 2, the main results
about strong and weak zero asymptotics are formulated, along with
some preliminary definitions and lemmas which are proved in
Section 3. In Section 4, a full set of orthogonality relations on
a single contour allows to pose a Riemann-Hilbert problem and to
apply the Deift and Zhou's steepest descent method (see
\cite{DeiftZhou}) to transform the Riemann-Hilbert problem in
order to obtain strong asymptotics of the solution. Finally, the
last section is devoted to the proofs of the main results.

\section{Main results}

\subsection{Basic definitions} \label{sec:basic_definitions}

Let us denote by $\C^+=\{z\in \C:\, \Im (z)>0 \}$ and $\C^-=\{z\in
\C:\, \Im (z)<0 \}$. For $A,B$ satisfying (\ref{region}), we
define the points
\begin{equation}\label{zeta}
    \zeta_{1,2}  =
    \frac{B^2 - A^2 \pm 4 \sqrt{(A+1)(B+1)(A+B+1)}}
    {(A+B+2)^2}.
\end{equation}
We will use the following convention: for $(A,B)$ such that
$A<-1<A+B$ (case C.2), $\zeta_1\in \C^+$ and $\zeta_2
=\overline{\zeta_1}$; for $-1<A<0<B$ (case C.3), we agree that
$-1<\zeta_1 <\zeta_2 <1$.

\begin{figure}[htb]
\centering \includegraphics[scale=0.8]{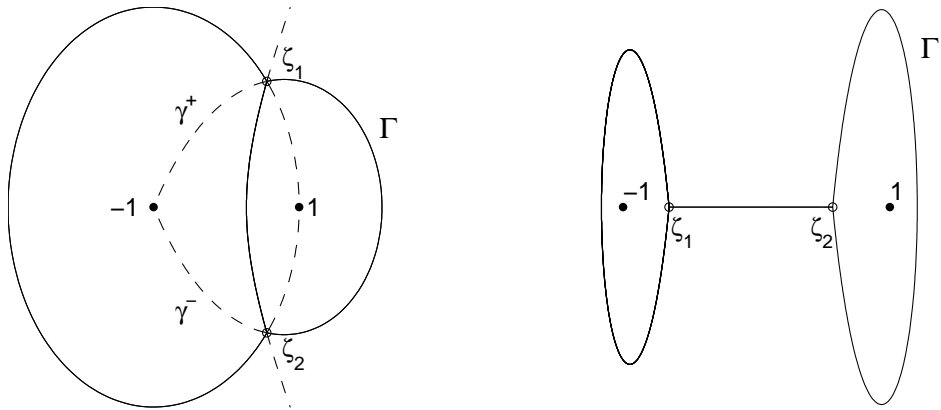}
\caption{Typical structure of the set $\mathcal N$ for cases C.2
(left) and C.3. Dashed lines on the left are orthogonal critical
trajectories.} \label{fig:setN}
\end{figure}

With these $\zeta_{1,2}$, consider the set
\begin{equation}\label{trajectories}
\mathcal{N} \isdef  \left\{ z \in \mathbb{C}  :\,\Re
\int_{\zeta_1}^z \frac{((t-\zeta_1)(t-\zeta_2))^{1/2}}{t^2-1} \,
dt = 0 \right\},
\end{equation}
where we continue the integrand analytically along the path of
integration. Obviously, the set does not depend on the branch of
the square root. In fact, it coincides with the union of the
critical trajectories of the quadratic differential
\begin{equation} \label{quaddiff}
    -\frac{(z-\zeta_1)(z-\zeta_2)}{(z^2-1)^2} \, dz^2,
\end{equation}
or more precisely, with their projection on $\C$. Taking into
account the local structure of trajectories of quadratic
differentials (see e.g. \cite{Pommerenke} or \cite{strebel}), we
can prove the following (see Figure \ref{fig:setN}):
\begin{lemma} \label{lemma:trajectories}
If parameters $A$, $B$ satisfy condition \eqref{region}, then for
$\zeta_{1,2}$ defined in \eqref{zeta} the quadratic differential
is regular. In other words, all its critical trajectories are
finite and have the following global structure (see Figure
\ref{fig:setN}):
\begin{itemize}
    \item For $(A,B)$ such that $A<-1<A+B$ (case C.2), $\NN$
    consists of three arcs which connect $\zeta_{1,2}$ and intersect the real line in exactly one
point, in such a way that each of the intervals $(-\infty,-1)$,
$(-1,1)$, $(1,\infty)$ is cut by one of these arcs.
    \item For $(A,B)$ such that $-1<A<0<B$ (case C.3), $\NN$
    consists of three arcs; one of them is the real interval
    $[\zeta_1,\zeta_2]$ and the other two are Jordan contours, passing
    through $z=\zeta_1$ (respect., $z=\zeta_2$) and enclosing $z=-1$
    (respect., $z=1$).
\end{itemize}
\end{lemma}

Now we define some relevant curves. We denote by $\Gamma$ the
rightmost curve from $\NN$. For case C.2, $\Gamma$ consists of an
arc connecting $\zeta_{1,2}$ and crossing once the interval
$(1,+\infty)$, and for case C.3, it is a closed contour passing
through $z=\zeta_2$ and surrounding $z=1$. For case C.2 we also
consider  the orthogonal trajectories $\NN ^\perp$ (defined by
replacing $\Im $ in (\ref{trajectories}) by $\Re $). As in Lemma
\ref{lemma:trajectories}, it is easy to prove that their global
structure is as appears in Figure \ref{fig:setN}, left (dashed
lines). We denote by $\gamma^+$ the arc of $\NN ^\perp$ joining
$\zeta_1$ and $-1$, and $\gamma^-=\overline{\gamma^+}$.

\begin{figure}[htb]
\centering \includegraphics[scale=0.8]{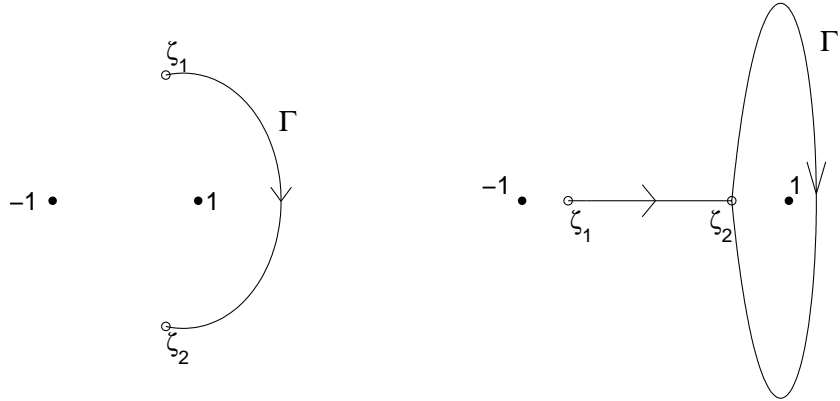}
\caption{Contours $\Sigma$ for cases C.2 (left) and C.3.}
\label{fig:contourSigma}
\end{figure}

Finally, we define the set $\Sigma$ as the smallest connected
subset of $\NN$ containing $\zeta_{1,2}$ and $\Gamma$. Namely,
\begin{equation}\label{Sigma}
\Sigma \isdef \begin{cases} \Gamma,\,& \quad \text{if }  A<-1<A+B, \text{ (case C.2)}\,;\\
\Gamma \cup [\zeta_1,\zeta_2],& \quad \text{if } -1<A<0<B, \text{
(case C.3)} \,.
\end{cases}
\end{equation}
As we see, in case C.2 the set $\Sigma$ is made of one critical
trajectory of the quadratic differential \eqref{quaddiff}, while
in case C.3 it is made of two. In both cases $\Sigma$ is oriented
from $\zeta_1$ to $\zeta_2$, and, in case C.3, clockwise, in such
a way that $(1, +\infty)$ is cut from the upper to the lower
half-plane (see Figure \ref{fig:contourSigma}).
For any function $f$ analytic and single-valued in a neighborhood
of $\Sigma$, this selection of the orientation induces two
boundary values of $f$ on $\Sigma$ that we denote by $f_+$ and
$f_-$ depending if we approach $\Sigma$ from the left or from the
right, respectively. On the sequel, we shall make use of the
concept of the polynomial convex hull of $\Sigma$, which is
denoted by $\Pc(\Sigma)$. In case C.2, $\Pc(\Sigma)=\Gamma $, so
that $\Int (\Pc(\Sigma)) =\emptyset$, where by $\Int(e)$ we denote
the set of inner points of $e$. Analogously, in case C.3,
$\Pc(\Sigma)$ is the union of $\Sigma$ and of the closure of the
bounded component of its complement, given by $\Int
(\Pc(\Sigma))$.

Next we define some functions that will play a role in what
follows. We denote
$$
R(z) \isdef \sqrt{(z-\zeta_1)(z-\zeta_2)}.
$$
It is a multi-valued and analytic function in $\C$, and we select
its single-valued branch in a plane cut from $\zeta_1$ to
$\zeta_2$ by imposing that
$$
\lim_{z\to \infty} \frac{R(z)}{z}=1\,.
$$
This allows us to define the (a priori complex) measure:
\begin{equation}\label{measure}
    d\mu(z)=\frac{A+B+2}{2\pi i}\;\frac{R_{+}(z)}{1-z^2}\;dz,\quad
    z\in \Sigma\,,
\end{equation}
with $\Sigma $ defined in \eqref{Sigma} and oriented as explained.
By Lemma \ref{lemma measure} below, $\mu$ is a unit positive
measure on $\Sigma$. Using Cauchy's Theorem it is easy to find an
analytic expression for its Cauchy transform:
\begin{equation}\label{hatmuC2}
\widehat{\mu}(z) \isdef \int \frac{d
\mu(t)}{z-t}=\frac{A+B+2}{2}\,
    \frac{R(z)}{z^2-1}-\frac{A/2}{z-1}-\frac{B/2}{z+1}\,, \quad z
    \in \C \setminus \Pc(\Sigma)\,;
\end{equation}
additionally, in case C.3,
\begin{equation}\label{hatmuC3}
\widehat{\mu}(z)=-\frac{A+B+2}{2}\,
    \frac{R(z)}{z^2-1}-\frac{A/2}{z-1}-\frac{B/2}{z+1}\,, \quad z
    \in \Int( \Pc(\Sigma))\,.
\end{equation}

Now, let us a define in $\C\setminus \Sigma$ a function which
plays a key role in the description of the strong asymptotics of
Jacobi polynomials,
\begin{equation}\label{G}
G(z) \isdef \exp\left(\int^{z}_{\zeta_2}\widehat{\mu}(t)\,
dt\right)\,.
\end{equation}
We normalize $G$ by imposing that
\begin{equation}\label{z0}
\lim_{z\to \zeta_2}G(z) = 1\,,
\end{equation}
where the limit is taken with $z$ approaching $\zeta _2$ from
$\C\setminus \Gamma$ (in case C.2) or from $\C^+\setminus \Sigma$
(in case C.3). Observe that since $\widehat \mu $ is the Cauchy
transform of a unit measure on $\Sigma$, function $G$ is analytic
and single-valued in $\C \setminus \Sigma$ in both cases
considered. Taking into account \eqref{hatmuC2}--\eqref{hatmuC3},
we see that there exists
\begin{equation}\label{defKappa}
\kappa \isdef \lim_{z\to \infty}\frac{G(z)}{z}\,.
\end{equation}
In addition, let
\begin{equation}\label{w}
w(z)= w(z;A,B) \isdef c(z-1)^{A/2} (z+1)^{B/2}=\exp\left(
\int^{z}_{\zeta_2}\left( \frac{A/2}{t-1}+\frac{B/2}{t+1}\right)\,
dt \right)\,,
\end{equation}
which is a multi-valued analytic function in $\mathbb{C}\setminus
\{\pm 1\}$. In what follows, we fix its single-valued analytic
branch in $\mathbb{C}\setminus (-\infty,1]$, by choosing the
constant $c$ (or the path of integration) in such a way that
\begin{equation}\label{normalizationW}
\lim_{ z\to \zeta_2} w(z) = 1 \,,
\end{equation}
where again the limit is taken with $z$ approaching $\zeta _2$
from $\C\setminus \Gamma$ (in case C.2) or from $\C^+\setminus
\Sigma$ (in case C.3).

The last ingredient for the asymptotics is given by the functions
\begin{equation}\label{defNN}
N_{11}(z)\isdef \frac{a(z)+a(z)^{-1}}{2} \quad \text{and} \quad
N_{12}(z) \isdef \frac{a(z)-a(z)^{-1}}{2i}
\end{equation}
(this notation is chosen because they will be entries of a certain
matrix $N$, see \eqref{model}), where
\begin{equation}\label{a}
a(z) \isdef
 \left(\frac{z-\zeta_{2}}{z-\zeta_{1}}\right)^{\frac{1}{4}}
\end{equation}
is defined in $\mathbb{C}\setminus \Gamma$ for $A<-1<A+B$ (case
C.2), and in $\mathbb{C}\setminus [\zeta_1,\zeta_2]$ for
$-1<A<0<B$ (case C.3). We select the branch of $a$ imposing the
normalization condition
$$
\lim_{z\to \infty} a(z)=1\,.
$$
Then, $N_{11}(z)\to 1$ and $N_{12}(z)\to 0$ as $z\to \infty$.

\subsection{Strong asymptotics}

First, we consider the strong asymptotics for Jacobi polynomials
with varying parameters satisfying \eqref{AB}--\eqref{region} with
$z$ away from $\Sigma$.
\begin{theorem} \label{theoremstrongOutside}
Let $(A,B)$ satisfy \eqref{region}. Then, for  $n\rightarrow
\infty$, the monic Jacobi polynomials $p_n =\widehat{P}_{n}
^{(An,Bn)}$ have the following asymptotic behavior:
\begin{equation}\label{AsympOutside}
p_{n}(z)=\left(\frac{G(z)}{\kappa}\right)^n\,N_{11}(z)\,\left(1+O\left(\frac{1}{n}\right)\right)\,,
\end{equation}
locally uniformly in $\C \setminus \Pc(\Sigma)$, where constant
$\kappa$ was defined in \eqref{defKappa}.

Furthermore, in the bounded component of $\C \setminus \Sigma$,
\begin{equation}\label{AsympC3Gamma}
\begin{split}
 p_n
(z)= &
\frac{1}{\kappa^n}\,\left(\left(G(z)w^2(z)\right)^{-n}N_{11}(z)
        \left(1 + O\left(\frac{1}{n}\right)\right)\right. \\ &
        \left. +2ie^{-An \pi i} \sin (A\pi n)\, G^n(z) N_{12}(z)
        \left(1 + O\left(\frac{1}{n}\right)\right) \right)\,.
\end{split}
\end{equation}
\end{theorem}
In particular, this theorem shows that zeros of $P_n^{(A n, B n)}$
do not accumulate in $\C \setminus \Pc(\Sigma)$.

Next, we describe the asymptotics on $ \Sigma $, but away from the
branch points $\zeta_{1,2}$:
\begin{theorem} \label{theoremstrongInside}
Let $(A,B)$ satisfy \eqref{region}. Then, for $n\rightarrow
\infty$, the monic Jacobi polynomials $p_n =\widehat{P}_{n}
^{(An,Bn)}$ have the following asymptotic behavior for $z$ away
from $\zeta_{1,2}$ :
\begin{enumerate}
\item[\rm (a)] In case C.2, on the ``$\pm$''-side of $\Gamma$
\begin{equation} \label{AsympC2Gamma}
\begin{split}
p_n (z) = & \frac{1}{\kappa^n}\,\left(G^n(z) N_{11}(z)
        \left(1 + O\left(\frac{1}{n}\right)\right)\right. \\
        & \left. \pm\,\left(G(z)w^2(z)\right)^{-n}N_{12}(z)
        \left(1 + O\left(\frac{1}{n}\right)\right) \right)\,.
\end{split}
\end{equation}

\item[\rm (b)] In case C.3, on the ``$-$''-side of $\Gamma$
formula \eqref{AsympC3Gamma} is still valid, while on the
``$+$''-side of $\Gamma$,
\begin{equation}\label{AsympC3GammaBis}
\begin{split}
 p_n
(z)= & \frac{1}{\kappa^n}\,\left( G^n(z) N_{11}(z)
        \left(1 + O\left(\frac{1}{n}\right)\right)\right. \\ &
        \left. +2ie^{-An \pi i} \sin (A\pi n)\, \left(G(z)w^2(z)\right)^{-n} N_{12}(z)
        \left(1 + O\left(\frac{1}{n}\right)\right) \right)\,.
\end{split}
\end{equation}

\item[\rm (c)] In case C.3, on the ``$+$''-side of
 $(\zeta_1,\zeta_2)$,
\begin{equation} \label{AsympC3int}
\begin{split}
p_n (z) = & \frac{1}{\kappa^n}\,\left(G^n(z) N_{11}(z)
        \left(1 + O\left(\frac{1}{n}\right)\right) + \left(G(z)w^2(z)\right)^{-n}N_{12}(z)
        \left(1 + O\left(\frac{1}{n}\right)\right) \right)\,,
\end{split}
\end{equation}
while on the ``$-$''-side of  $(\zeta_1,\zeta_2)$,
\begin{equation} \label{AsympC3intDown}
\begin{split}
p_n (z) = & \frac{1}{\kappa^n}\,\left(G^n(z) N_{11}(z)
        \left(1 + O\left(\frac{1}{n}\right)\right) \right. \\ & - \left.
        e^{-2A \pi i n} \, \left(G(z)w^2(z)\right)^{-n}N_{12}(z)
        \left(1 + O\left(\frac{1}{n}\right)\right) \right)\,.
\end{split}
\end{equation}
\end{enumerate}
\end{theorem}

\begin{remark} \label{remarkstrongC2}
All these asymptotic expressions match on the boundaries of the
overlapping domains and on the respective regions. For instance,
it will be shown that in a small neighborhood of every point of
$\Sigma$ (distinct from $\zeta _{1,2}$), $|G(z) w(z)|>1$ for
$z\notin \Sigma$. Hence, the first term in
\eqref{AsympC2Gamma}--\eqref{AsympC3intDown} is dominant, and away
from $\Sigma$ they reduce to \eqref{AsympOutside}.

Furthermore, in the case C.3, if $An$ are not exponentially close
to integers (in the sense that will be made more precise below),
the second term in \eqref{AsympC3Gamma} dominates, and we may
write
\begin{equation*} \label{AsympC3Gammared}
p_n(z)=\,\left(\frac{G(z)}{\kappa}\right)^n\,\,2i e^{-An\pi i}\,
\sin (A\pi n)N_{12}(z)
        \left(1 + O\left(\frac{1}{n}\right)\right)\,.
\end{equation*}
\end{remark}

\medskip

The asymptotic formulas above are no longer valid close to the
branch points $\zeta_{1, 2}$. As it usually happens in a
neighborhood of the ``soft ends'' of the support of an equilibrium
measure, asymptotics is described in terms of the Airy function
$\Ai(z)$ and its derivative. We give explicit formulas only for
the rightmost (according to the orientation of $\Sigma$) branch
point $\zeta _2$, which is in a certain sense, the ``interesting''
one. Clearly, the analysis at the other point is similar.

In order to formulate our result in a more compact form it is
convenient to introduce the function
\begin{equation}\label{phiC2}
    \phi(z)\isdef \frac{A+B+2}{2}\int_{\zeta_{2}}^{z}
    \frac{R(t)}{1-t^{2}}\, dt\,.
\end{equation}
Consider case C.2. Given a sufficiently small $\varepsilon >0$,
and a neighborhood $\Delta_\varepsilon (\zeta _2)\isdef \{z\in
C:\, |z-\zeta _2|<\varepsilon  \}$, it is a single-valued analytic
function in $\Delta_\varepsilon (\zeta _2) \setminus \Gamma$.
Furthermore, taking into account the local behavior of $R$,
function
\begin{equation}\label{f}
    f(z)\isdef \frac{3}{2}\, \left( \phi (z)\right)^{2/3}
\end{equation}
can be extended as single-valued to the whole $\Delta_\varepsilon
(\zeta _2)$. Here the $2/3$rd power is chosen such that $f(z)>0$
along $\gamma^-$.

\begin{theorem} \label{theoremstrongC2endpoint}
Let $(A,B)$ such that $A<-1<A+B$ (case C.2). Then, there exists $
\varepsilon >0$ such that if $|z-\zeta_2|< \varepsilon $, we have
that the monic Jacobi polynomials $p_n =\widehat{P}_{n}
^{(An,Bn)}$ satisfy:
\begin{equation*}\label{strongnearC2}
\begin{split}
p_{n}(z)=\frac{ \sqrt{\pi}}{\kappa^n w^n(z)}\,&  \left(
\frac{n^{1/6}\, f^{1/4}(z)}{a(z)}\, \Ai (n^{2/3} f(z)) \, \left(
1+O\left( \frac{1}{n}\right) \right) \right. \\ %
 & \left. - \frac{a(z)}{n^{1/6}\, f^{1/4}(z)}\, \Ai' (n^{2/3} f(z)) \, \left(
1+O\left( \frac{1}{n}\right) \right) \right)\,,
\end{split}
\end{equation*}
where $a(z)$ is defined in \eqref{a}, and we take $f^{1/4}(z)>0$
along $\gamma^-$.
\end{theorem}

\begin{remark}
Obviously, the asymptotic behavior near $\zeta _1$ in this case is
completely symmetric to $\zeta_2$ with respect to $\R$.
\end{remark}

Consider case C.3. For a sufficiently small $0<\varepsilon
<1-\zeta _2$, $\phi$ is single-valued and analytic in
$\Delta_\varepsilon (\zeta _2) \setminus (\zeta _1, \zeta _2)$.
Function $f$, defined again by formula \eqref{f}, can be extended
as a single-valued function to the whole $\Delta_\varepsilon
(\zeta _2)$, with the $2/3$rd power chosen such that $f(z)>0$
along $(\zeta _2, 1)$.

\begin{theorem} \label{theoremstrongC3endpoint}
Let $(A,B)$ such that $-1<A<0<B$ (case C.3). Then, there exists $
\varepsilon >0$ such that if $|z-\zeta_2|< \varepsilon $, we have
that the monic Jacobi polynomials $p_n =\widehat{P}_{n}
^{(An,Bn)}$ satisfy:
\begin{equation} \label{localCaseC3}
\begin{split}
p_{n}(z)=\frac{\sqrt{\pi}}{\kappa^n w^n(z)}\,&  \left(
\frac{n^{1/6}\, f^{1/4}(z)}{a(z)}\, \mathcal A (n^{2/3} f(z)) \,
\left(
1+O\left( \frac{1}{n}\right) \right) \right. \\ %
 & \left. - \frac{a(z)}{n^{1/6}\, f^{1/4}(z)}\, \mathcal A' (n^{2/3} f(z)) \, \left(
1+O\left( \frac{1}{n}\right) \right) \right)\,,
\end{split}
\end{equation}
where $a(z)$ is defined in \eqref{a},
$$
\mathcal A(t)= \mathcal{A}(t;A,n) \isdef e^{-A\pi i n} \, \Ai(t) +
2ie^{ \pi i/ 3 } \, \sin(A\pi n) \, \Ai\left( e^{ 4 \pi i/ 3 } \,
t \right)\,,
$$
and we take $f^{1/4}(z)>0$ along $(\zeta _2, 1)$.
\end{theorem}

\begin{remark}
Formulas stated in Theorems
\ref{theoremstrongOutside}--\ref{theoremstrongC3endpoint} are
locally uniformly continuous both on the $z$ and $(A,B)$ planes,
which allows to extend them to the general case of $\{\alpha _n\}$
and $\{ \beta _n\}$ satisfying \eqref{AB}--\eqref{region}.
\end{remark}

\subsection{Weak zero asymptotics}%
\label{subsec:Weak zero asymptotics}

As a corollary of the asymptotic formulas stated in the previous
section we can obtain the distribution of the zeros of the
sequence of polynomials $P_n^{(\alpha_n, \beta_n)}$, where $\{
\alpha_n\}$ and $\{ \beta _n\}$ satisfy
\eqref{AB}--\eqref{region}. By ``weak zero asymptotics'' we
understand here the limit (in the weak-* sense) of the normalized
zero counting measures associated with $P_n^{(\alpha_n,
\beta_n)}$.

\begin{figure}[htb]
\centering \includegraphics[scale=0.6]{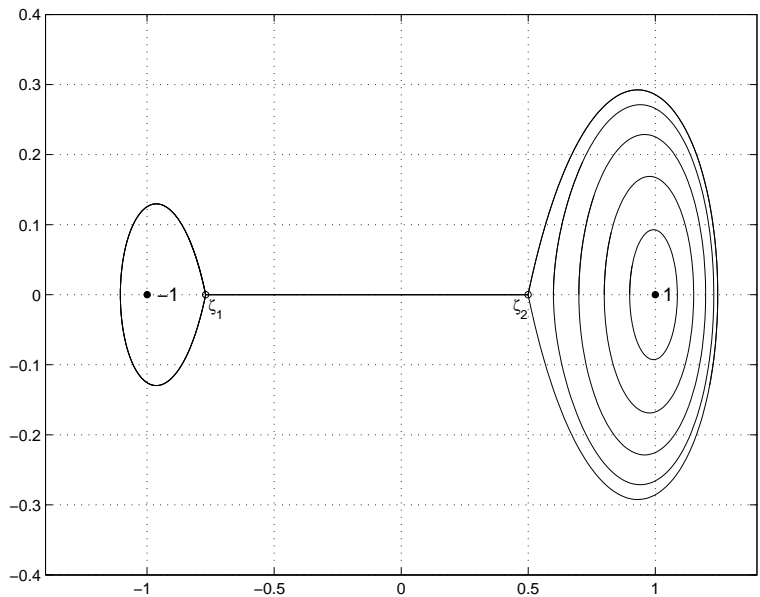}
\caption{Some trajectories of the quadratic differential
(\ref{quaddiff}), or equivalently, some level sets $\Gamma_r$, for
the values $A = -0.8$ and $B=0.5$. } \label{fig:trajectories}
\end{figure}

The measure $\mu$ introduced in \eqref{measure} will be all we
need for the description of the asymptotic behavior of the zeros
in the case C.2. However, region C.3 comprises the pathological
cases given by \eqref{integer 1}. By continuity, we may expect
here a variety of limit behaviors. In fact, in order to
characterize completely the weak zero asymptotics of Jacobi
polynomials in the case C.3  we need to use a 1-parametric family
of measures including (\ref{measure}). Namely, in case C.3, we
must consider the sets
\begin{equation*}
\NN_{r}\, =\NN^{(A,B)}_{r} \isdef \left\{z\in \mathbb{C}:\, |G(z)
w(z)|=e^{r/2}\right\}\,= \left\{ z \in \mathbb{C}  :\,\Re
\int_{\zeta_2}^z \frac{R(t)}{t^2-1} \, dt = \frac{r}{A+B+2}
\right\},
\end{equation*}
for $r\geq 0$. They also consist of trajectories of the quadratic
differential (\ref{quaddiff}), and $\mathcal N_0=\mathcal N$. Now,
we define $\Gamma_r$ as the rightmost curve in $\NN_{r}$ or, what
is the same, the part of $\NN_{r}$ which is entirely contained in
the half-plane $\{z\in \mathbb{C}:\, \Re z\geq \zeta_2\}$. It is
easy to check that for $r>0$ the level curve $\Gamma_{r}$ is a
closed contour inside $\Gamma =\Gamma_{0}$ surrounding the point
$z=1$ (see Figure \ref{fig:trajectories}).

For each $r\in [0,\infty)$ we define the absolutely continuous
measure
\begin{equation}\label{measure level}
d\mu_{r}(z)\isdef \frac{A+B+2}{2\pi
i}\;\frac{R_{+}(z)}{1-z^{2}}\;dz,\quad
    z\in (\zeta_{1},\zeta_{2})\cup \Gamma_{r},
\end{equation}
and, for $r=\infty$, the measure
\begin{equation}\label{measure Dirac}
d\mu_{\infty}(z)\isdef
 -A\delta_{1} + \frac{A+B+2}{2\pi}\frac{\sqrt{(z-\zeta_{1})(\zeta_{2}-z)}}{1-z^{2}}\,
 \chi_{[\zeta_{1},\zeta_{2}]}\,dz,
\end{equation}
where $\chi_{[\zeta_{1},\zeta_{2}]}$ is the characteristic
function of the interval $[\zeta_{1},\zeta_{2}]$, and $\delta_1$
is the Dirac delta (unit mass point) at $z=1$.

\begin{lemma} \label{lemma measure}
If $(A, B)$ satisfy \eqref{region}, then for $r \geq 0$ measures
$\mu_r $ in \eqref{measure level}--\eqref{measure Dirac} (and, in
particular, measure $\mu$ in (\ref{measure})) are unit positive
measures. Moreover, for $(A,B)$ satisfying $-1<A<0<B$ (case C.3),
we have that  for $0\leq r\leq +\infty$,
\begin{equation}\label{measC3}
    \mu_r ([\zeta_1,\zeta_2]) = 1+A,\quad \mu_r(\Gamma_r) = -A\,,
\end{equation}
where we consider $\Gamma_0=\Gamma$ and $\Gamma_\infty=\{1\}$.
\end{lemma}

\medskip

Now we are ready to state the weak zero asymptotics for the Jacobi
polynomials $P_n^{(\alpha _ n ,  \beta _ n)}$. In the case C.3,
when  $-1<A<0<B$, we make an additional assumption: the sequence
of parameters $\alpha_n$ satisfies that the following limit
exists:
\begin{equation}\label{integerdist}
\lim_{n \to \infty} \left[ \dist(\alpha_n, \mathbb
Z)\right]^{1/n}\,=\, e^{-r}\,, \quad 0\leq r\leq +\infty.
\end{equation}
Then, it holds:
\begin{theorem} \label{theoremweak}
Consider a sequence of Jacobi polynomials
$P_n^{(\alpha_n,\beta_n)}$, $n\in \mathbb{N}$, such that sequences
$\{\alpha_n\}$, $\{\beta_n\}$ satisfy (\ref{AB}) and
(\ref{region}). Then:
\begin{enumerate}

\item[\rm (i)] If $(A, B)$ satisfy condition C.2, then the zeros
of $P_n^{(\alpha_n,\beta_n)}$, $n\in \mathbb{N}$, accumulate on
the arc $\Gamma$, and measure $\mu$ in (\ref{measure}) is the
weak* limit of the corresponding normalized zero counting
measures.

\item[\rm (ii)] If $(A, B)$ satisfy condition C.3, and
(\ref{integerdist}) holds for some $0 \leq r\leq +\infty$, then
the zeros of Jacobi polynomials $P_n^{(\alpha_n,\beta_n)}$, $n\in
\mathbb{N}$, accumulate on $[\zeta_1,\zeta_2]\cup \Gamma_r$, and
measure $\mu_r$ defined above is the weak* limit of the normalized
zero counting measures.

\end{enumerate}
\end{theorem}

As we said before, part (i) of Theorem \ref{theoremweak},
corresponding to case C.2, was established in \cite{MR2002d:33017}
for parameters $\alpha_n,\beta_n$ varying according to (\ref{AB})
but with $A,B<-1$, which is a region of the $(A,B)$-plane
equivalent to $A<-1<A+B$ by means of transformations
(\ref{symmetry 1})--(\ref{symmetry 2}).

In the case C.3 the situation when $r=0$, that is,
$$
\lim_{n \to \infty} \left[ \dist(\alpha_n, \mathbb
Z)\right]^{1/n}\,=\,1,
$$
is generic, because it takes place when parameters  do not
approach the integers exponentially fast. When the limit in
(\ref{integerdist}) is smaller than one, i.e. $r>0$, curve
$\Gamma$ is replaced by a level curve $\Gamma_r$, strictly
contained inside the bounded component of the complement to
$\Gamma$, and the support of the limit measure becomes
disconnected. Finally, when parameters $\alpha_n , n\in
\mathbb{N}$, tend to integers faster than exponentially, the limit
measure has a discrete part consisting of a Dirac mass at $z=1$.
We could have anticipated this phenomenon observing the
coalescence of zeros given by \eqref{integer 1}.

\begin{figure}
\centering
\begin{tabular}{c}
\mbox{\includegraphics[scale=0.7]{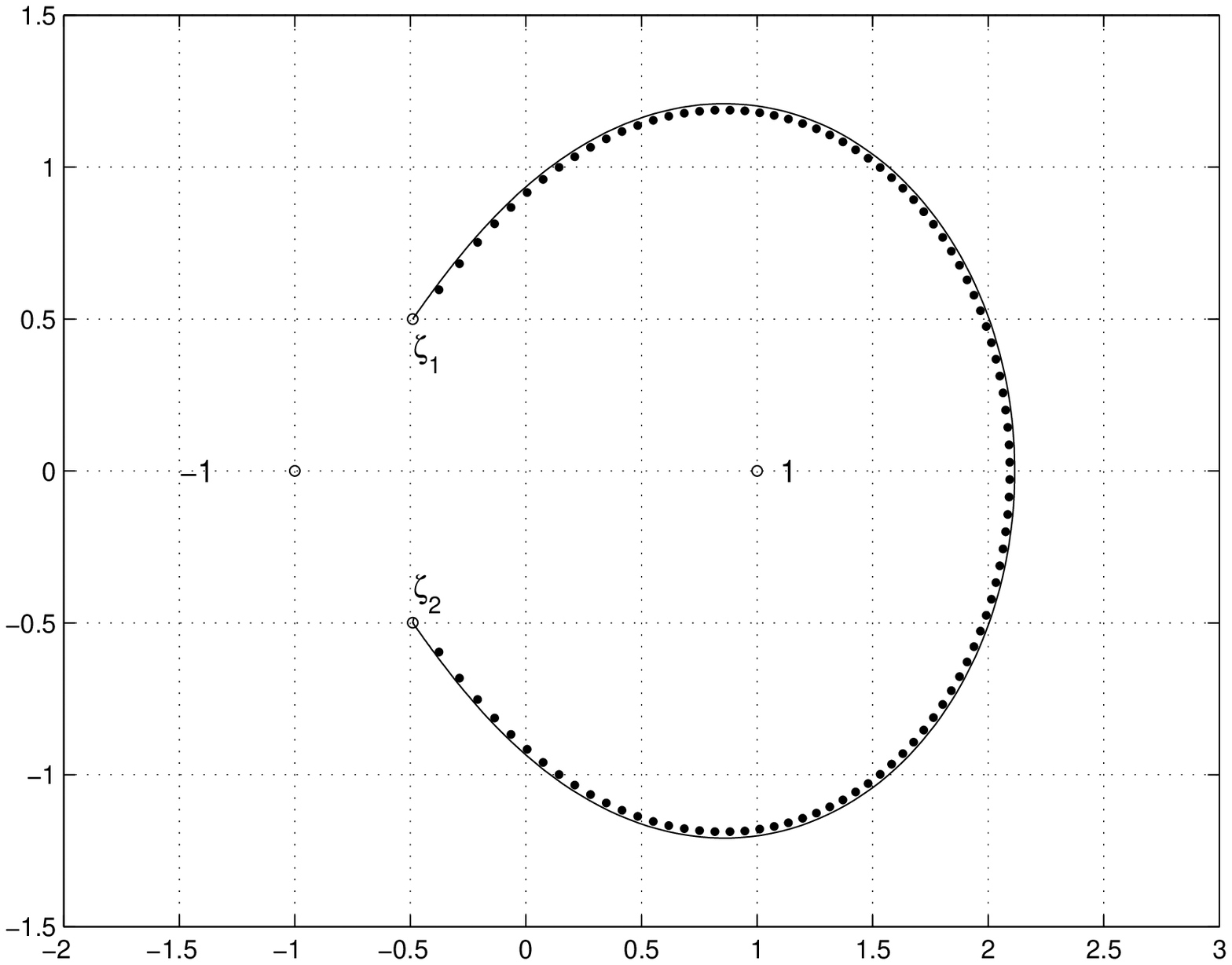}} \\
\mbox{\includegraphics[scale=0.7]{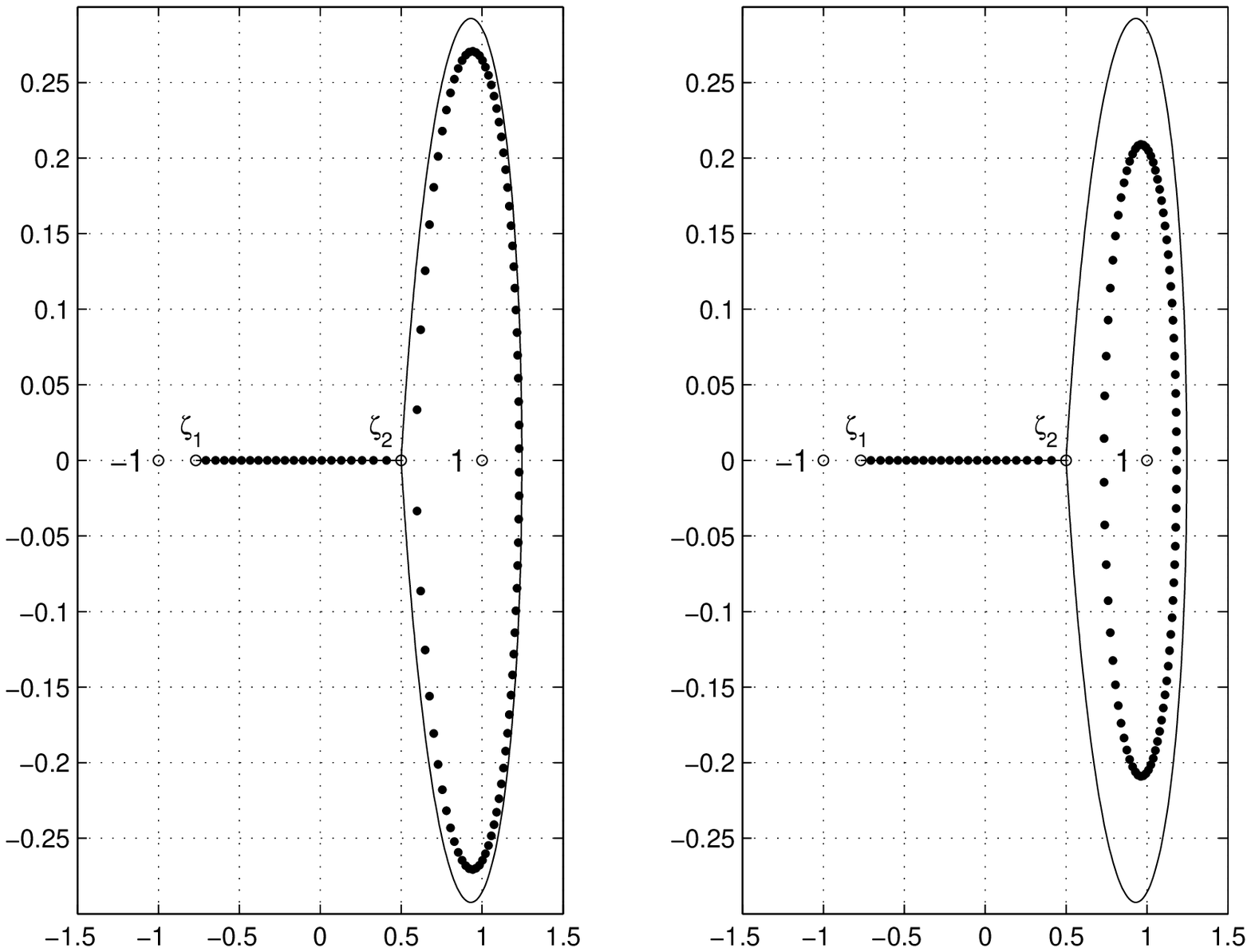}}
\end{tabular}
\caption{Up, case C.2: zeros of $P_{100}^{(-110 - 10^{-5}, 50 -
10^{-5})}$, along with the curve $\Gamma$, corresponding to
$A=-1.1$ and $B=0.5$. Down, case C.3: zeros of $P_{100}^{(-80 -
10^{-5}, 50 -  10^{-5})}$ (left) and $P_{100}^{(-80 - 10^{-15}, 50
-  10^{-5})}$ (right), along with the set $\Sigma$, corresponding
to $A=-0.8$ and $B=0.5$.} \label{fig:C3}
\end{figure}

Examples of zeros of Jacobi polynomials for Cases C.2 and C.3 are
represented in Figure \ref{fig:C3}.

\begin{remark}
The case C.3 deals with situations when real zeros arise. When
$-n<\alpha_n <0$ and $\beta_n >-1$, $P_n^{(\alpha,\beta)}$ satisfy
a quasi-orthogonality relation (see Theorem 6.1 in \cite{Etna})
which ensures the existence of, at least, $n-[-\alpha_n]$ zeros in
$(-1,1)$. This lower bound of the number of zeros in $(-1,1)$ is
exact, according to the so-called Hilbert-Klein formulas
\cite[Theorem 6.72]{szego:1975}. Since $\lim_{n\rightarrow
\infty}\frac{n-[-\alpha_n]}{n}\,=\,1+A$, looking at \eqref{measC3}
we see that the mass of the part of the asymptotic measure of
zeros supported on $[\zeta_1,\zeta_2] \subset (-1,1)$ agrees with
the limit of the ratio of zeros placed in $(-1,1)$, given by the
Hilbert-Klein formulas.
\end{remark}

\begin{remark}
At this point, it is natural to ask about what happens when $A=-1$
and $B>0$, which is a transition case between C.2 and C.3. By
\eqref{symmetry 1}--\eqref{symmetry 2}, it describes also the
situation when $(A,B)$ belongs to any of the straight lines
$A=-1$, $B=-1$ and $A+B=-1$, outside of the square $(A,B)\in
[-1,0]\times [-1,0]$. Roughly speaking, in this case the endpoints
$\zeta_{1, 2}$ are confluent in a single point, say $\zeta$, and
zeros approach a simple closed contour emanating from this point
(critical trajectory) and surrounding $z=1$, or other closed
trajectories strictly contained in the interior of this critical
trajectory. In \cite{MR2002d:33017}, the equation of this critical
trajectory is conjectured. This conjecture is proved in
\cite{Driver99} and \cite{Duren01} for the particular case where
$\alpha_n =-n-1$ and $\beta_n =kn+1$, $k$ being a fixed positive
real number. This and the other transitions between cases C.1 --
C.5 deserve a separated treatment.
\end{remark}

\section{Proof of the auxiliary lemmas}

\begin{varproof} \textbf{of Lemma \ref{lemma:trajectories}.}
This proof is based upon the local structure of the trajectories
of quadratic differentials (see \cite{Pommerenke} or
\cite{strebel}). We restrict our attention to the case where
$-1<A<0<B$ (case C.3). The proof of the other case is similar (see
also the proof of Lemma 2.1 in \cite{ArnoAndrei}).

First, we see that for $-1<A<0<B$ the quadratic differential
\eqref{quaddiff} possesses two simple zeros at $\zeta_{1, 2}$.
Thus, we know that three critical trajectories emanate from
$\zeta_{1, 2}$ at equal angles. Moreover, the segment $[\zeta _1,
\zeta _2]\subset \NN$, which is straightforward to verify by
definition of $\NN$.

On the other hand, \eqref{quaddiff} has double poles at $z=\pm 1$
and $z=\infty$, in such a way that if we consider the rational
function $Q(z)=-(z-\zeta_1)(z-\zeta_2) (z^2-1)^{-2}$, we have that
the residues of $\sqrt{Q}$ at these points are purely imaginary.
Therefore, we conclude that near these double poles the
trajectories are simple closed contours.

Now, the symmetry of $Q$ with respect to the real axis, along with
the facts that the trajectories cannot tend to infinity and
\eqref{quaddiff} has no other singular point, allows one to ensure
that the other critical trajectories are two closed contours
emanating from $\zeta_1$ and $\zeta_2$. The fact that a closed
trajectory needs to surround a singular point implies that these
closed trajectories intersect the real axis in two points, one of
them in $(1,+\infty)$ and the other in $(-\infty,-1)$.

\end{varproof}

\begin{varproof}\textbf{of Lemma \ref{lemma measure}.}
Taking into account the definition of $R(z)$  and (\ref{zeta}), it
is easy to see that:
\begin{equation}\label{R(1)}
R(1)= \begin{cases} \dfrac{2A}{A+B+2}<0, & \quad \text{if }
   A<-1<A+B \text{ (case C.2)}\,;\\
-\dfrac{2A}{A+B+2}>0,& \quad \text{if }  -1<A<0<B
 \text{ (case C.3)}\,.
\end{cases}
\end{equation}
In the same way,
\begin{equation}\label{R(-1)}
R(-1)\,=\,-\frac{2B}{A+B+2}<0, \quad \text{if }  A<0<B \quad
\text{and} \quad A+B>-1 \,.
\end{equation}
Thus, the definition of $\Sigma$ in (\ref{Sigma}) yields that
(\ref{measure}) is real-valued on $\Sigma$ and does not change
sign on each of its components. The same remains valid when
$-1<A<0<B$ and $0<r\leq +\infty$ for $\mu_r$ on its support
$\Sigma_r =\Gamma_r \cup [\zeta_1,\zeta_2]$.

Moreover, for $(A,B)$ such that $-1<A<0<B$, taking into account
the definition of $\Sigma$ and $\Sigma_r$, the residue theorem and
(\ref{R(1)})--(\ref{R(-1)}), we have for $0\leq r<+\infty$:
\begin{equation*}\label{intGamma}
\begin{split}
\mu_r(\Gamma_r) &= \int_{\Gamma_r}d\mu_r (t) =
-(A+B+2)\res_{z=1}\left( \frac{R(z)}{z^2-1}\right) =
(A+B+2)\frac{R(1)}{2} = -A\, ;
\end{split}
\end{equation*}
clearly, also $\mu_{\infty}(\Gamma_{\infty}) =
\mu_{\infty}(\{1\})\,=\,-A$.

On the other hand, for $0\leq r\leq +\infty$,
\begin{equation*}\label{intzeta}
\begin{split}
\mu_r([\zeta_1,\zeta_2])\,& =\, \int_{\zeta_1}^{\zeta_2}
d\mu_r(t) \\
& =\,\frac{A+B+2}{2}\left[\res_{z=1}\left(
\frac{R(z)}{z^2-1}\right) +\res_{z=-1}\left(
\frac{R(z)}{z^2-1}\right)+\res_{z=\infty}\left(
\frac{R(z)}{z^2-1}\right)\right] \\
&
=\,\frac{A+B+2}{2}\left(1+\frac{A}{A+B+2}-\frac{B}{A+B+2}\right)\,=\,1+A\,
,
\end{split}
\end{equation*}
and, therefore,

\begin{equation*}\label{intSigma}
\mu_r(\Sigma_r)\,=\,\int_{\Sigma_r}d\mu_r(t) =
\int_{\Gamma_r}d\mu_r(t) + \int_{\zeta_1}^{\zeta_2}
d\mu_r(t)\,=\,1\,.
\end{equation*}
Analogously, for $(A,B)$ such that $A<-1<A+B$, it is easy to see
that
\begin{equation*}\label{intC2}
\mu(\Gamma)\,=\,\int_{\Gamma}d\mu(t) = 1\,,
\end{equation*}
and it settles the proof.
\end{varproof}

\section{Riemann-Hilbert analysis}

\subsection{Orthogonality and the Riemann-Hilbert problem}
As it was mentioned in the Introduction, the key fact for the
asymptotic analysis is a full system of orthogonality relations
satisfied by the Jacobi polynomials on simple contours, which
allows to pose a matrix Riemann-Hilbert problem (RHP) and apply
the Deift-Zhou steepest descent method.

The following result was established in \cite[Theorem 5.1]{Etna}:
\begin{lemma} \label{lemma:etna}
Let $\CC$ be a Jordan arc connecting $z=-1+0i$ with $z=-1-0i$ and
surrounding $z=1$ once.  If $\beta >0$, then we have:
\begin{equation}\label{orthNonHerm}
    \int_{\CC} t^k\, P_n^{(\alpha , \beta )}(t)w^2(t;\alpha,\beta) \,
    dt \begin{cases}
    =0\,, & k<n,\\
    \neq 0 \,, & k=n,
    \end{cases}
\end{equation}
where $w(\cdot;\alpha,\beta)$ has been introduced in \eqref{w}.
\end{lemma}

From the seminal work of Fokas, Its and Kitaev \cite{Fokas92} (see
also \cite{Deift99}) it is known that the orthogonality
\eqref{orthNonHerm} can be characterized in terms of the following
Riemann-Hilbert problem: find a matrix valued function $Y :\,
\mathbb{C}\backslash\CC \rightarrow\mathbb{C}^{2 \times 2}$
satisfying the conditions below:
\begin{enumerate}
 \item[(RH1.1)] $Y$ is analytic in $\C \setminus \CC$.
 \item[(RH1.2)] $Y$ has continuous boundary values on $\CC$, denoted by $Y_+$
and $Y_-$, such that
\[ Y_+(z) = Y_-(z)\begin{pmatrix}
  1 & w^2(z; \alpha, \beta) \\
  0 & 1 \\
\end{pmatrix}\,,
 \quad \text{for } z\in \CC. \]
 \item[(RH1.3)] As
$z\to\infty$,
\[ Y(z) = \left(I + O\left(\frac{1}{z}\right)\right)
\begin{pmatrix}
  z^{n} & 0 \\
  0 &  z^{-n}\\
\end{pmatrix} \,.\]
\item[(RH1.4)] $Y$ is bounded in a neighborhood of $z=-1$.
\end{enumerate}
\begin{proposition}[\cite{Etna}] \label{solRHY} The unique solution of the
Riemann-Hilbert problem (RH1.1)--(RH1.4) is given by
\begin{equation*} \label{formulaY}
Y(z)=\begin{pmatrix}
  p_{n}(z) & \frac{1}{2\pi i} \int_{\CC} \frac{p_{n}(t)w^2(t; \alpha, \beta)}
        {t-z}\, dt \\
  q_{n-1}(z) &  \frac{1}{2\pi i} \int_{\CC} \frac{q_{n-1}(t)
        w^2(t; \alpha, \beta)}{t-z}\, dt\\
\end{pmatrix}\,,
\end{equation*}
where $p_{n}(z)=\widehat P_n^{(\alpha , \beta )}(z)$ is the monic
Jacobi polynomial, and $q_{n-1}(z)=b_{n-1}P_{n-1}^{(\alpha , \beta
)}(z)$, for some suitable non-zero constant $b_{n-1}$.
\end{proposition}

Let $(A, B)$ be a pair satisfying \eqref{region}. Then for every
$n\in \N$, monic polynomials $\widehat P_n^{(An , Bn)}$ satisfy
the conditions of Lemma \ref{lemma:etna}. Hence, for any $\CC$ as
described in Lemma \ref{lemma:etna}, polynomials $\widehat
P_n^{(An , Bn)}$ can be characterized as the $(1,1)$ entry of the
unique matrix $Y$ solving the Riemann-Hilbert problem
(RH1.1)--(RH1.4) with
\begin{equation*}\label{alpha_beta}
    \alpha = A n\,, \quad \beta =B n \,.
\end{equation*}

Taking advantage of the freedom in the selection of the contour
$\CC$ in \eqref{orthNonHerm} we will choose different
configurations for both cases C.2 and C.3. This choice is mainly
suggested by the numerical evidence on the actual location of
zeros. For describing the appropriate $\CC$ we will use the
contours defined in subsection \ref{sec:basic_definitions},
corresponding to $A$ and $B$ fixed.

\begin{figure}[htb]
\centering \includegraphics[scale=0.8]{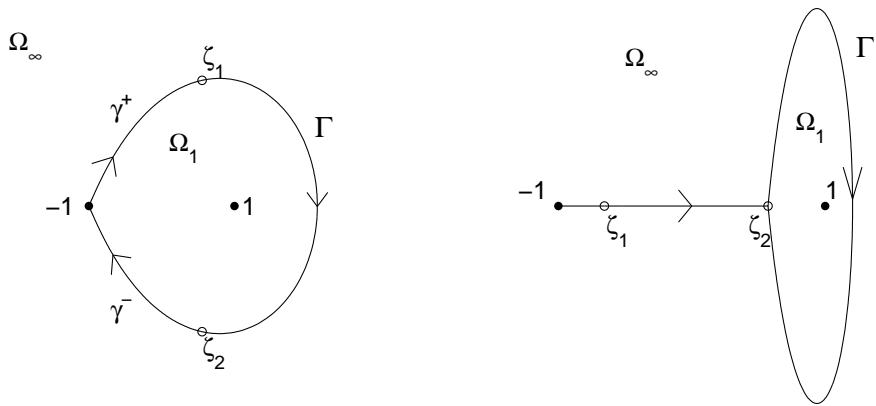}
\caption{Contours $\CC$ for the Riemann-Hilbert problem for $Y$;
cases C.2 (left) and C.3.} \label{fig:contourInt}
\end{figure}

In the case C.2, when $A<-1< A+B$, we will make the contour $\CC$
in Lemma \ref{lemma:etna} coincide with $\gamma^- \cup \Gamma \cup
\gamma^+$, oriented clockwise (Figure \ref{fig:contourInt}, left).
Hence, we are interested in the Riemann-Hilbert problem
(RH1.1)--(RH1.4) with $\CC=\gamma^- \cup \Gamma \cup \gamma^+$,
$\alpha = A n$ and $ \beta =B n$.

However, in case C.3, when $-1<A<0<B$, it is convenient to make
part of the original contour $\CC$ coalesce along the interval
$[-1, \zeta_2]$ (traversed twice in opposite directions), and the
rest go along the arc $\Gamma$. This deformation creates a new
contour, which we denote again by $\CC$, and we choose its
orientation as in Figure \ref{fig:contourInt}, right.
Consequently, it yields a new RHP, still characterizing the
polynomials $p_n$. With respect to the problem (RH1.1)--(RH1.4),
we have to modify only the jump matrix on $(-1, \zeta _2)$ as a
result of the coalescence of two original sub-arcs of $\CC$: its
$(1,2)$ entry becomes $w_+^{2n}-w_-^{2n}$ on $(-1, \zeta _2)$. The
new Riemann-Hilbert problem is: find a matrix valued function $Y
\equiv Y^{(A,B)}:\mathbb{C}\setminus ([-1, \zeta _2] \cup
\Gamma)\rightarrow\mathbb{C}^{2 \times 2}$ such that the following
conditions hold:
\begin{enumerate}
\item[(RH2.1)] $Y$ is analytic on $\mathbb{C}\setminus ([-1, \zeta
_2] \cup \Gamma)$.
 \item[(RH2.2)] $Y$ has continuous boundary
values on $\CC \setminus \{-1, \zeta _2\}$, denoted by $Y_+$ and
$Y_-$, such that $Y_+(z) = Y_-(z)J_{Y}(z)$, where
\[ J_{Y}(z) = \begin{cases}
    \begin{pmatrix}
  1 &  w^{2n}(z)\\
  0 & 1 \\
    \end{pmatrix}, & z\in \Gamma; \\
    \begin{pmatrix}
  1 &  d_n  \,  w^{2n}_+(z)
  \\
  0 & 1 \\
    \end{pmatrix}, & z\in (-1,\zeta_{2})\,;  \end{cases}\]
with
\begin{equation}\label{defDn}
d_n \isdef 1-e^{-2An \pi i}=2 i e^{-An \pi i}\,\sin(An\pi)\, ,
\end{equation}
and with $w(z)=w(z;A,B)$ defined in \eqref{w}.

 \item[(RH2.3)] As $z\to\infty$,
\[ Y(z) = \left(I +O\left(\frac{1}{z}\right)\right)
\begin{pmatrix}
  z^{n} & 0 \\
  0 &  z^{-n}\\
\end{pmatrix}\,.\]
\item[(RH2.4)] $Y$ is bounded in a neighborhood of $z=-1$ and
$z=\zeta _2$.
\end{enumerate}

In both cases $\mathbb{C}\setminus \CC$ has two connected
components, one containing $z=1$ and the other containing
infinity; we denote these components by $\Omega_1$ and
$\Omega_\infty$, respectively (see Figure \ref{fig:contourInt}).

The steepest descent analysis, that we are going to carry out
next, introduces new contours which are unions of a finite number
of curves and arcs on $\C$. In order to simplify notation we will
call all the end points and points of self-intersection of such
curves \emph{singular} points, and the rest will be \emph{regular}
points of the contour. Hence, we could rephrase (RH1.4) and
(RH2.4) saying that $Y$ is bounded in a neighborhood of all
singular points of $\CC$.

\subsection{First transformation $Y \mapsto U$}

\begin{figure}[htb]
\centering \includegraphics[scale=0.8]{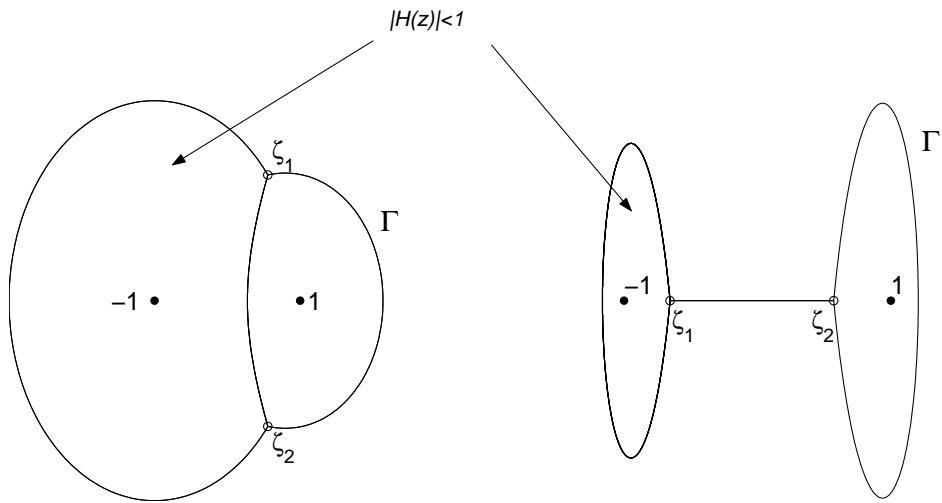}
\caption{Regions where $|H(z)|<1$ for cases C.2 (left) and C.3.}
\label{fig:regionsCases23}
\end{figure}

In order to shorten notation, we use the Pauli matrix
\begin{equation*}\label{Pauli}
\sigma_{3}\isdef \begin{pmatrix}
    1 & 0 \\
    0 & -1 \\
\end{pmatrix},
\end{equation*}
and denote $x^{\sigma_3} = \begin{pmatrix}
    x & 0 \\
    0 & x^{-1} \\
\end{pmatrix}$.
Also for the sake of brevity, it is convenient to introduce the
function
\begin{equation}\label{H}
H(z) \isdef G(z) w(z)\,,
\end{equation}
analytic and single-valued in $\C \setminus (\Sigma \cup (-\infty,
1])$. Using \eqref{hatmuC2}--\eqref{normalizationW}, we
immediately get that
\begin{itemize}
\item For the case C.2,
 \begin{equation}\label{GtimesWcaseC2}
H(z)
=\exp\left(\dfrac{A+B+2}{2}\int_{\zeta_2}^z\frac{R(t)}{t^2-1}dt\right)\,,
\quad \text{for}\, z\in \mathbb{C}\setminus \left( \Gamma \cup
(-\infty, 1]\right) \,.
\end{equation}
\item For the case C.3,
 \begin{equation}\label{GtimesW}
H(z) = \begin{cases} \displaystyle
\exp\left(\dfrac{A+B+2}{2}\int_{\zeta_2+i0}^z\frac{R(t)}{t^2-1}dt\right),\,
& \quad \text{for}\, z\in \mathbb{C}\setminus \left( \Pc (\Sigma)\cup (-\infty, 1]\right) \,;\\
 \displaystyle \exp\left(-\frac{A+B+2}{2}\int_{\zeta_2}^z\frac{R(t)}{t^2-1}dt\right),\,
& \quad \text{for}\, z\in \Int(\Pc (\Sigma))\cap \C^+ \,,\\
 \displaystyle e^{-\pi i A}\, \exp\left(-\frac{A+B+2}{2}\int_{\zeta_2}^z\frac{R(t)}{t^2-1}dt\right),\,
& \quad \text{for}\, z\in \Int(\Pc (\Sigma))\cap \C^- \,.
\end{cases}
\end{equation}
Observe that the same convention as in \eqref{z0} for the path of
integration applies:
$$
\lim_{\C^+ \setminus \Sigma \ni z\to \zeta_2} H(z) = 1 \,.
$$
Furthermore, taking into account \eqref{measC3}, in the case C.3,
\begin{equation}\label{valueAtzeta}
\lim_{\C^- \setminus \Pc(\Sigma) \ni z\to \zeta_2} H(z) = e^{ -\pi
i \mu(\Gamma) }=e^{A \pi i }\,.
\end{equation}
\end{itemize}
In both cases C.2 and C.3, the sets of trajectories $\mathcal{N}$
and $\mathcal{N}_r$, $(r\geq 0)$, introduced in subsection
\ref{subsec:Weak zero asymptotics}, may be characterized by the
conditions $|H(z)|=1$ and $|H(z)|=e^{r/2}$, respectively. Figure
\ref{fig:regionsCases23} shows  also the domains where $|H(z)|<1$.

Function $H$ has continuous boundary values at regular points of
$\Sigma \cup (-\infty, 1]$, which satisfy:
\begin{itemize}
\item In case C.2:
$$
H_+(z)=\begin{cases}
    H_-^{-1}(z), & z\in \Gamma\,, \\
    e^{\pi i A}\, H_-(z), & z \in (-1,1)\,, \\
    e^{\pi i (A+B)}\, H_-(z), & z \in (-\infty,-1)\,.
    \end{cases}
$$
\item In case C.3:
\begin{equation}\label{boundaryValuesH}
    H_+(z)=\begin{cases}
    H_-^{-1}(z), & z\in \Gamma\,, \\
    e^{- \pi i A}\, H_-^{-1}(z), & z\in (\zeta _1, \zeta _2)\,, \\
    e^{\pi i A}\, H_-(z), & z \in (-1,1)\setminus [\zeta _1, \zeta _2]\,, \\
    e^{\pi i (A+B)}\, H_-(z), & z \in (-\infty,-1)\,.
    \end{cases}
\end{equation}
\end{itemize}
This allows us to express the boundary values of $G$ at the
regular points of $\Sigma$ in terms of $H$:
\begin{equation}\label{bvG1}
\frac{G_+(z)}{G_-(z)}=H_+^2(z)\,, \quad \text{and} \quad G_+(z)\,
G_-(z)=\frac{1}{w^2_+(z)}\,.
\end{equation}

Now we are ready to introduce the first transformation of the RHP,
with the aim to normalize it at infinity. For $d_n$ in
\eqref{defDn}, let us fix any value of $d_n^{1/2}$, and define
\begin{equation}\label{first transf}
    U(z)\,=\,\begin{cases}
    \kappa^{n\sigma_3}Y(z)G(z)^{-n\sigma_3}, & \quad \text{in case
    C.2}\,; \\
     d_n ^{-\sigma_3/2}\kappa^{n\sigma_3}Y(z)G(z)^{-n\sigma_3} d_n ^{\sigma_3/2},
    & \quad \text{in case C.3}\,,
    \end{cases}
\end{equation}
with $\kappa $ given by \eqref{defKappa}. Obviously, matrix $U$
solves now a new Riemann-Hilbert problem. Taking into account
\eqref{bvG1} we can state it as:
\begin{enumerate}
\item[(RH3.1)] $U$ is analytic on $\mathbb{C}\setminus \CC$.
\item[(RH3.2)] $U$ has continuous boundary values at the regular
points of $\CC$, denoted by $U_+$ and $U_-$, such that $U_+(z) =
U_-(z)J_{U}(z)$, where, for $(A,B)$ in case C.2,
\begin{equation}\label{JU C2}
   J_U=\begin{cases}
    \begin{pmatrix}
  H_+^{-2n}(z) & 1 \\
  0 & H_+^{2n}(z)
  \end{pmatrix}\,, & z \in \Gamma\,, \\
  \begin{pmatrix}
  1 & H^{2n}(z)   \\
  0 & 1
  \end{pmatrix}\,, & z \in \gamma^+ \cup \gamma^- \,, \\
  \end{cases}
\end{equation}
and for $(A,B)$ in case C.3,
\begin{equation}\label{JU C3}
J_{U}(z) = \begin{cases}
    \begin{pmatrix}
  H_+^{-2n}(z) &  d_{n}^{-1}\\
  0 & H_+^{2n}(z)  \\
    \end{pmatrix}, & z\in \Gamma \,; \\
    \begin{pmatrix}
  H_+^{-2n}(z)  &  1\\
  0 &   H_+^{2n}(z)  \\
    \end{pmatrix}, & z\in (\zeta_{1},\zeta_{2});\\
    \begin{pmatrix}
    1 & H_+^{2n}(z)
    \\
    0 & 1 \\
    \end{pmatrix}, & z\in (-1,\zeta_{1}). \\
\end{cases}
\end{equation}

\item[(RH3.3)] As $z\to\infty$,

\[ U(z) =  I +O\left(\frac{1}{z}\right) \,.\]

\item[(RH3.4)] Matrix $U$ is bounded in a neighborhood of the
singular points of $\CC$.
\end{enumerate}

\subsection{Second transformation $U \mapsto T$}

By \eqref{JU C2}--\eqref{JU C3}, the jump matrix $J_U$ has
oscillatory diagonal entries on $\Sigma$, along with exponentially
decaying (as $n\to \infty$) off-diagonal entries elsewhere and
away from $\zeta_{1, 2}$ (see Figure \ref{fig:regionsCases23}).
The aim of the next step is to transform the jump matrices with
oscillatory diagonal entries into matrices asymptotically close to
the identity matrix or to matrices with constant jumps. To this
end, we take advantage of an appropriate factorization of $J_U$
and ``open the lenses'' around contours $\CC$.

In case C.2, we use the following factorization of the jump matrix
for $z\in \Gamma$ (where we have taken into account that
$H_+=1/H_-$ on $\Gamma$):
\begin{equation}\label{factor C2}
J_U(z)= \begin{pmatrix}
  1 & 0 \\
  H_-^{-2n}(z) & 1
  \end{pmatrix} \, \begin{pmatrix}
  0 & 1 \\
  -1 & 0
  \end{pmatrix} \, \begin{pmatrix}
  1 & 0 \\
  H_+^{-2n}(z) & 1
  \end{pmatrix}\,.
\end{equation}

\begin{figure}[htb]
\centering \includegraphics[scale=0.85]{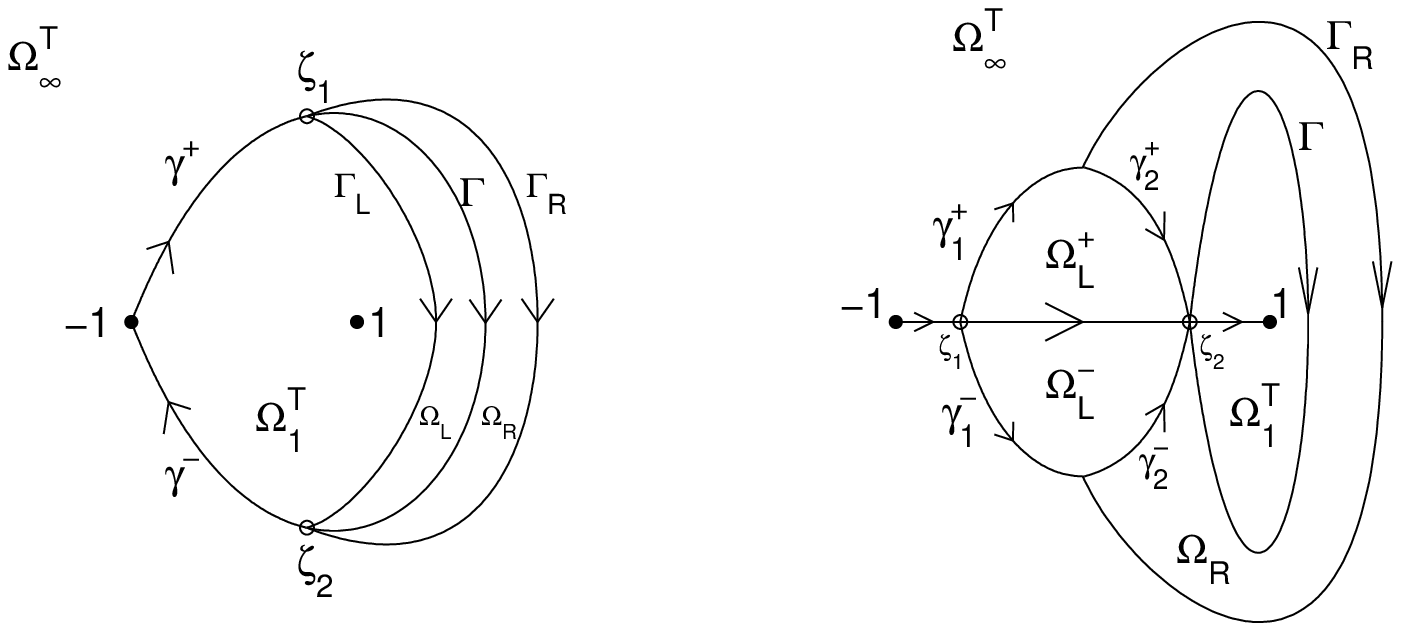}
\caption{Contours $\CC^T$ for $T$.} \label{fig:FirstLensesCases23}
\end{figure}

Thus, the problem of the oscillatory diagonal entries of the jump
matrix for $z \in \Gamma$ may be solved by opening the lenses
around $\Gamma$ as it is shown in Figure
\ref{fig:FirstLensesCases23} (left). The new contours $\Gamma_L$
and $\Gamma_R$ are also oriented from $\zeta_1$ to $\zeta_2$, and
this gives us two new bounded regions, $\Omega_R$ and $\Omega_L$,
as well as modified domains $\Omega_1 ^{T}\isdef \Omega_1
\setminus \overline{\Omega_L} $ and $\Omega_\infty^{T}\isdef
\Omega_\infty \setminus \overline{\Omega_R}$; we also denote
$\CC^T \isdef  \CC \cup \Gamma_L \cup \Gamma_R$, with the
orientation shown in Figure \ref{fig:FirstLensesCases23}, left.

Hence, taking into account (\ref{factor C2}), we define in case
C.2 a new matrix-valued function $T:\mathbb{C}\setminus\CC^{T}
\longrightarrow \mathbb{C}^{2\times 2}$ by
    \begin{equation}\label{change_var2}
T(z)\,=U(z)\cdot \,\begin{cases} I, & \text{for } z \in \Omega_1
^T \cup
\Omega_\infty ^T \,, \\
\begin{pmatrix}
  1 & 0 \\
  H^{-2n}(z) & 1
  \end{pmatrix}\,, & \text{for } z \in \Omega_L\,, \\
U(z)\, \begin{pmatrix}
  1 & 0 \\
  -H^{-2n}(z) & 1
  \end{pmatrix}\,,&  \text{for } z \in \Omega_R\,.
\end{cases}
\end{equation}
It solves the following RHP:
\begin{enumerate}
    \item[(RH4.1)] $T$ is analytic for $z \in \C \setminus \CC^T$;
    \item[(RH4.2)] $T(z)$ possesses continuous boundary values
    at regular points of $\CC^T $, $T_+$ and $T_-$, related by the following
    jump conditions:
    $$
   T_+(z)=T_-(z) J_T(z)\,, \quad z \in \CC^T\,,
   $$
   where the jump matrix $J_T$ is
      $$
   J_T (z)=\begin{cases}
   \begin{pmatrix}
  0 & 1  \\
  -1 & 0
  \end{pmatrix}\,, & z \in \Gamma\,, \\
  \begin{pmatrix}
  1 & H^{2n}(z)  \\
  0 & 1
  \end{pmatrix}\,, & z \in \gamma^+ \cup \gamma^- \,, \\
    \begin{pmatrix}
    1 & 0 \\
    H^{-2n}(z) & 1
  \end{pmatrix}\,, & z \in \Gamma_L \cup \Gamma_R\,.
   \end{cases}
   $$

    \item[(RH4.3)] $T(z)$ has the following behavior at
    infinity:
    $$
  T(z)= I + O(1/z) \quad \text{as } z \to \infty\,.
    $$
    \item[(RH4.4)] $T(z)$ is bounded in a neighborhood of the singular points of $\CC^T$.
\end{enumerate}

In principle, we could take advantage of factorization
\eqref{factor C2} also in the case C.3. However, the geometry here
is more complicated; this procedure would eventually yield a
constant jump on whole $\Sigma$, which has now two components. In
order to give a unified treatment to both cases in the next step,
we use now a different factorization for $J_U$:
\begin{equation*}\label{factor C3}
J_U(z)=
\begin{cases}
\begin{pmatrix} 0 & d_n^{-1 } \\ -d_n  &
H_-^{-2n}(z) \end{pmatrix}\,
\begin{pmatrix} 1 & 0 \\
d_n    H_+^{-2n}(z) & 1 \end{pmatrix}\,, & z \in \Gamma\,, \\
\begin{pmatrix}
  1 & 0 \\
  e^{-2 A\pi ni}\, H_-^{-2n}(z) & 1
  \end{pmatrix} \, \begin{pmatrix}
  0 & 1 \\
  -1 & 0
  \end{pmatrix} \, \begin{pmatrix}
  1 & 0 \\
   H_+^{-2n}(z) & 1
  \end{pmatrix}\,, & z \in  (\zeta _1, \zeta _2)\,.
\end{cases}
\end{equation*}
These factorizations suggest to open lenses in the way shown in
Figure \ref{fig:FirstLensesCases23}, right, which yields the new
contour $\CC^{T}$, splitting $\C$ into domains $\Omega_1^T$,
$\Omega_{\infty}^{T}$, $\Omega_L^\pm$, and $\Omega_R$, as shown.
Now we define the matrix-valued function
$T:\mathbb{C}\setminus\CC^{T} \to \mathbb{C}^{2\times 2}$ in the
following way:
    \begin{equation}\label{T}
    T(z)=U(z)\cdot\,
    \begin{cases}
    \begin{pmatrix} 1 & 0 \\
    -d_n   H^{-2n}(z) & 1 \end{pmatrix}\,, & z \in \Omega_{R}\,, \\
    \begin{pmatrix} 0 & d_n^{-1 } \\ -d_n  &
    H^{-2n}(z) \end{pmatrix}\,, & z \in \Omega_{1}^T\,, \\
    \begin{pmatrix}
  1 & 0 \\
  -  H^{-2n}(z) & 1
  \end{pmatrix}\,, & z \in \Omega_{L}^+\,, \\
    \begin{pmatrix}
  1 & 0 \\
    e^{ -2A\pi ni}\, H^{-2n}(z) & 1
  \end{pmatrix}\,, & z \in \Omega_{L}^-\,,\\
    I, & z \in \Omega_{\infty}^{T} \,.
    \end{cases}
\end{equation}
Then $T(z)$ solves the problem (RH4.1)--(RH4.4), with $J_U$
replaced by
   \[ J_{T}(z) = \begin{cases}
    \begin{pmatrix}
  0 &  1 \\
  -1 & 0 \\
    \end{pmatrix}\,, & z\in (\zeta_1,\zeta_2)\,, \\
    \begin{pmatrix}
    1 & 0 \\
    d_{n}H^{-2n}(z) & 1 \\
    \end{pmatrix}\,, & z\in \Gamma_{R}\,, \\
    \begin{pmatrix}
    1 & 0 \\
    H^{-2n}(z) & 1 \\
    \end{pmatrix}\,, & z\in \gamma_1^+ \cup \gamma_2^-\,, \\
    \begin{pmatrix}
    1 & 0 \\
    e^{-2 A\pi ni}H^{-2n}(z) & 1 \\
    \end{pmatrix}\,, & z\in \gamma_1^- \cup \gamma_2^+\,, \\
    \begin{pmatrix}
    1 & H_+^{2n}(z) \\
    0 & 1 \\
    \end{pmatrix}\,, & z\in (-1,\zeta_1)\,, \\
    \begin{pmatrix}
    1 &  e^{2 A\pi ni}H_+^{-2n}(z)\\
    0 & 1 \\
    \end{pmatrix}\,, & z\in (\zeta_2,1) \,, \\
    I\,, & z\in \Gamma \,.
    \end{cases}\]

\subsection{Construction of the parametrices}

Now, we can see (cf.\ Figure \ref{fig:regionsCases23}) that we
have a single open arc joining the branch points ($\Gamma$ in case
C.2, and $(\zeta_1,\zeta_2)$ in case C.3) where the jump matrix
$J_T$ is constant, and at a positive distance from these arcs,
$J_T$ is asymptotically exponentially close to the identity
matrix. Hence, by ignoring the ``close-to-identity'' jumps and
condition (RH4.4) we are lead to the following problem: find an
analytic matrix-valued function $N(z)=I+O(1/z)$, $z\to \infty$,
and having the jump
$$
N_+(z)=N_-(z)\, \begin{pmatrix} 0 & 1 \\ -1 & 0
\end{pmatrix}
$$
on $\Gamma$ (in case C.2) or on $(\zeta_1,\zeta_2)$ in case C.3,
with the orientation ``from $\zeta _1$ to $\zeta_2$'' chosen. A
solution of this model RHP, which is not unique in general, is
(cf.\ \cite[Ch. 7]{Deift99}):
\begin{equation}\label{model}
    N(z)= \begin{pmatrix}
    \frac{a(z)+a(z)^{-1}}{2} & \frac{a(z)-a(z)^{-1}}{2i} \\
    -\frac{a(z)-a(z)^{-1}}{2i} & \frac{a(z)+a(z)^{-1}}{2}
    \end{pmatrix}\,,
\end{equation}
where $a$ has been defined in \eqref{a}; it satisfies
$$
N(z)=O(|z-\zeta_j|^{-1/4})\,, \quad z \to \zeta_j\,, \quad j =1,
2\,,
$$
showing that the singularities at $\zeta_j$ are
$L^{2}$-integrable. Observe that the $(1,1)$ and $(1,2)$ entries
of $N$ coincide with $N_{11}$ and $N_{12}$, introduced in
\eqref{defNN}.

We may expect $N$ to be close to $T$ away from $\zeta_1$ and
$\zeta_2$. However, in a neighborhood of the branch points the
ignored jumps are no longer close to identity, and a different
parametrix (model problem) is required. Now we look for two
matrices $P^{(j)}$, $j\in \{ 1, 2 \}$, which have the same jumps
as $T$ in a neighborhood of $z= \zeta _j$, and match $N$ on the
boundary of these neighborhoods.

The construction of these matrices is well described for instance
in \cite{Deift99}. Denote by $\Delta_{\varepsilon}(s)\isdef \{z\in
\mathbb{C}:\, |z-s|<\varepsilon\}$, where we take $\varepsilon >0$
sufficiently small. A local parametrix $P^{(j)}$ in
$\Delta_{\varepsilon}( \zeta _j)$, $j\in \{1, 2 \}$, solves the
RHP with the same jumps as $T$ there (see Figure
\ref{fig:localAnalysisCases2and3}):
\begin{figure}[htb]
\centering
\includegraphics[scale=0.7]{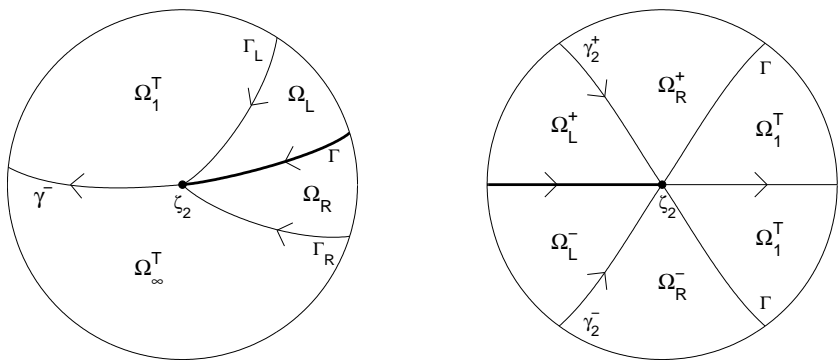}
\caption{Local analysis for cases C.2 (left) and C.3.}
\label{fig:localAnalysisCases2and3}
\end{figure}
\begin{enumerate}
    \item[(RH5.1)] $P^{(j)}$ is analytic for $z \in \Delta_{\varepsilon}( \zeta _j) \setminus \CC^T
    $, bounded and continuous in $\overline{\Delta_{\varepsilon}( \zeta _j)} \setminus \CC^T $;
    \item[(RH5.2)] $P^{(j)}(z)$ possesses continuous boundary values
    at regular points of $ \CC^T \cap
   \Delta_{\varepsilon}( \zeta _j)$, $P_+^{(j)}$ and $P_-^{(j)}$, related by the following
    jump conditions:
    $$
   P^{(j)}_+(z)=P^{(j)}_-(z) J_{P^{(j)}}(z),  \quad z \in \CC ^T \cap
   \Delta_{\varepsilon}( \zeta _j)\,.
   $$
   \item[(RH5.3)] there exists a constant $M >0$ such that for every
   $z \in \partial\Delta_{\varepsilon}( \zeta _j)\setminus \CC^{T}$,
    $$
  \| P^{(j)}(z) N^{-1}(z)- I \| \leq \frac{M}{n}\,.
    $$
\end{enumerate}
We describe the construction for $\zeta _2$; in order to simplify
notation we write $P$ instead of $P^{(2)}$ whenever it cannot lead
us into confusion. The jumps $J_P=J_{P^{(2)}}$ specified in
(RH5.2) are
\begin{itemize}
\item In case C.2:
   $$
   J_{P}(z)=\begin{cases}
   \begin{pmatrix}
  0 & 1  \\
  -1 & 0
  \end{pmatrix} \,, \quad z \in \Gamma \cap
  \Delta_{\varepsilon}(\zeta_2)\,, \\
\begin{pmatrix}
  1 & 0  \\
  H^{-2n}(z) & 1
  \end{pmatrix} \,, \quad z \in (\Gamma_L \cup \Gamma_R)\cap \Delta_{\varepsilon}(\zeta_2) \,,\\
\begin{pmatrix}
  1 & H^{2n}(z)  \\
  0 & 1
  \end{pmatrix} \,, \quad z \in \gamma^- \cap \Delta_{\varepsilon}(\zeta_2)\,.
   \end{cases}
   $$
\item In case C.3:
 \[ J_{P}(z) = \begin{cases}
    \begin{pmatrix}
  0 &  1 \\
  -1 & 0 \\
    \end{pmatrix}\,, & z\in (\zeta_2-\varepsilon ,\zeta_2)\,, \\
    \begin{pmatrix}
    1 & 0 \\
    e^{-2 A\pi ni}H^{-2n}(z) & 1 \\
    \end{pmatrix}\,, & z\in \gamma_2^+ \cap
    \Delta_{\varepsilon}(\zeta_2)\,,  \\
    \begin{pmatrix}
    1 & 0 \\
    H^{-2n}(z) & 1 \\
    \end{pmatrix}\,, & z\in  \gamma_2^- \cap
    \Delta_{\varepsilon}(\zeta_2)\,, \\
    \begin{pmatrix}
    1 &  e^{2 A\pi ni}H_+^{-2n}(z)\\
    0 & 1 \\
    \end{pmatrix}\,, & z\in (\zeta_2,\zeta_2+\varepsilon ) \,, \\
    I\,, & z\in \Gamma \,.
    \end{cases}\]
\end{itemize}

In order to solve the Riemann-Hilbert problems for $P$, let us
first make a simple change of functions yielding piecewise
constant jump matrices. For this purpose, we set for $z\in
    \overline{\Delta_{\varepsilon}(\zeta_2)}\setminus
    \CC^{T}$,
\begin{equation}\label{Rcommon}
 R(z)\isdef P(z) \cdot \begin{cases}
  e^{-n\phi(z) \sigma_3}, & \text{in case C.2,} \\
  e^{A \pi i n \sigma_3}\, e^{-n\phi(z) \sigma_3}, & \text{in case C.3,}
 \end{cases}
\end{equation}
where $\phi$ is the function introduced in \eqref{phiC2}. In order
to compute the new jumps we need to find how $\phi$ is related to
$H$. In the case C.2, by \eqref{GtimesWcaseC2}, $\exp(-\phi(z))=
H(z)$ for $z\in \Delta_{\varepsilon}(\zeta_2) \setminus \Gamma$.
In the case C.3, by continuity of $\phi$ in
$\Delta_{\varepsilon}(\zeta_2)\setminus (\zeta _1, \zeta _2)$,
\begin{equation}\label{relationPhiH}
\exp(-\phi(z))=\begin{cases} H(z)\,, & z\in \left( \Omega_L \cup
\Omega_R \right)\cap  \Delta_{\varepsilon}(\zeta_2) \cap \C^+\,,
\\%
H^{-1}(z)\,, & z\in  \Omega_1^T \cap \Delta_{\varepsilon}(\zeta_2)
\cap \C^+\,,
\\%
e^{-\pi i A}\, H^{-1}(z)\,, & z\in  \Omega_1^T \cap
\Delta_{\varepsilon}(\zeta_2) \cap \C^-\,,
\\%
e^{-\pi i A}\, H(z)\,, & z\in \left( \Omega_L \cup \Omega_R
\right)\cap \Delta_{\varepsilon}(\zeta_2) \cap \C^-\,,
\end{cases}
\end{equation}
and
$$
(\phi _+ + \phi _-)(z)=2\pi i A\,, \quad z \in (\zeta
_2-\varepsilon , \zeta _2)\,.
$$
Now we can compute the jump matrix for $R$: $J_{ R}=e^{n\phi_-(z)
\sigma_3}\, J_T \, e^{-n\phi_+(z) \sigma_3}$, namely:
 \begin{itemize}
\item in case C.2,
$$
J_{R}(z) = \begin{cases}
    \begin{pmatrix}
    0 & 1 \\
    -1 & 0 \\
    \end{pmatrix}, & z \in \Gamma \cap
    \Delta_{\varepsilon}(\zeta_2)\,, \\
    \begin{pmatrix}
    1 & 0 \\
    1 & 1 \\
    \end{pmatrix}, & z \in (\Gamma_L \cup \Gamma_R) \cap \Delta_{\varepsilon}(\zeta_2)\,, \\
    \begin{pmatrix}
    1 & 1 \\
    0 & 1 \\
    \end{pmatrix}, & z \in \gamma_-\cap \Delta_{\varepsilon}(\zeta_2)\,. \\
\end{cases}
$$
\item in case C.3,
$$
J_{ R}(z) = \begin{cases}
    \begin{pmatrix}
    0 & 1 \\
    -1 & 0 \\
    \end{pmatrix}, & z \in (\zeta_{2}-\varepsilon, \zeta_{2})\,,\\
    \begin{pmatrix}
    1 & 0 \\
     1 & 1 \\
    \end{pmatrix}, & z \in \gamma_{2}^{\pm}\cap
    \Delta_{\varepsilon}(\zeta_2)
    \,, \\
    \begin{pmatrix}
    1 & 1 \\
    0 & 1 \\
    \end{pmatrix}, & z \in (\zeta_{2}, \zeta_{2}+\varepsilon)\,.
     \\
    \end{cases}
$$
\end{itemize}
Observe now that we have essentially the same local problem in
both cases; in this way, we have reduced the RHP to the one
studied in \cite{Deiftetal1} (see also \cite[Ch.\ 7]{Deift99}),
and we can write its solution explicitly.

For $\varepsilon >0$ small enough, function
$$
    f(z)= \frac{3}{2}\, \left( \phi (z)\right)^{2/3},
$$
defined in \eqref{f}, is a conformal mapping from the neighborhood
of the branch point onto a neighborhood of 0. In case C.2,
$\Gamma$ and $\gamma^-$ are mapped onto the negative and positive
real axis, respectively, and in case C.3 it happens to $(\zeta _1,
\zeta _2)$ and $(\zeta _2,1)$. Also we may deform the other curves
($\Gamma_L$ and $\Gamma_R$ in case C.2, and $\gamma_2^{\pm}$ in case
C.3) in such a way that the points on their image by $f$ close to
the branch point have the argument $\pm 2\pi /3$.

Then the problem for $R$ is solved by
$$
R(z)=\Psi \left( n^{2/3} f(z)\right)\,,
$$
where $\Psi$ is built out of the Airy function $\Ai$ (see e.g.
\cite{Abramowitz}) and its derivative $\Ai'$ as follows:
\begin{equation}\label{Airy}
\Psi(t)\,=\,\begin{cases}
    \begin{pmatrix}
    \Ai(t) & \Ai(\omega^{2}t) \\
    \Ai'(t) & \omega^{2}\Ai'(\omega^{2}t) \\
    \end{pmatrix} e^{-\frac{\pi i}{6} \sigma_{3}}, & 0<\arg t<2\pi /3 ;\\
    \begin{pmatrix}
    \Ai(t) & \Ai(\omega^{2}t) \\
    \Ai'(t) & \omega^{2}\Ai'(\omega^{2}t) \\
    \end{pmatrix} e^{-\frac{\pi i}{6}\sigma_3} \begin{pmatrix}
    1 & 0 \\
    -1 & 1 \\
    \end{pmatrix}, & 2\pi /3 <\arg t<\pi;\\
    \begin{pmatrix}
    \Ai(t) & -\omega^{2} \Ai(\omega t) \\
    \Ai'(t) & -\Ai'(\omega t) \\
    \end{pmatrix} e^{-\frac{\pi i}{6}\sigma_3} \begin{pmatrix}
    1 & 0 \\
    1 & 1 \\
    \end{pmatrix}, & -\pi <\arg t<-2\pi /3;\\
    \begin{pmatrix}
    \Ai(t) & -\omega^{2} \Ai(\omega t) \\
    \Ai'(t) & -\Ai'(\omega t) \\
    \end{pmatrix} e^{-\frac{\pi i}{6}\sigma_3}, & -2\pi /3 <\arg
    t<0,
    \end{cases}
\end{equation}
and $\omega\isdef e^{2\pi i/3}$. Finally, matrix $P$ solving
(RH5.1)--(RH5.3) is
\begin{equation}\label{Pcommon}
 P(z)= E(z)\, R(z) \cdot \begin{cases}
  e^{n\phi(z) \sigma_3}, & \text{in case C.2,} \\
  e^{-A \pi i n \sigma_3}\, e^{n\phi(z) \sigma_3}, & \text{in case C.3,}
 \end{cases}
\end{equation}
where the analytic matrix function $E$ is
\begin{equation}\label{En}
    E(z) \isdef \sqrt{\pi}e^{\frac{\pi i}{6}}\,
    \begin{pmatrix}
    1 & -1 \\
    -i & -i \\
    \end{pmatrix}\,
    \left( \frac{n^{1/6} f(z)^{1/4}}{a(z)} \right)^{\sigma_3}.
\end{equation}

\subsection{Final transformation $T \mapsto S$}

Now we may use matrix valued functions $N$ and $P^{(j)}$ for the
final transformation. Recalling the definition of the contour
$\CC^T$, define the matrix-valued function $S$:
\begin{equation}\label{V}
    S(z) \isdef \begin{cases}
    T(z)N(z)^{-1},\, & z \in \mathbb{C} \setminus (\CC^{T} \cup
    \overline{\Delta_{\varepsilon}(\zeta_1)} \cup
    \overline{\Delta_{\varepsilon}(\zeta_2)}) \\
    T(z)\left(P^{(j)}(z)\right)^{-1},\, & z \in \Delta_{\varepsilon}(\zeta_j)\,, \quad j=1, 2 \,.
    \end{cases}
\end{equation}

\begin{figure}[htb]
\centering \includegraphics[scale=0.85]{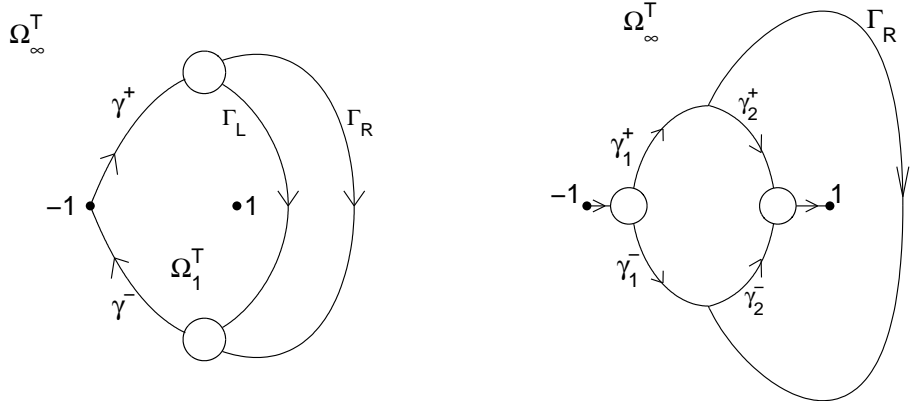}
\caption{Contours $\CC^S$ for $S$ in cases C.2 (left) and C.3.}
\label{fig:FinalLensesCases23}
\end{figure}
It is immediate to check that $S$ is analytic in
$\mathbb{C}\setminus \CC^{S}$, where $\CC^{S}$ is the contour
shown in Figure \ref{fig:FinalLensesCases23}. Moreover, $S:\,
\mathbb{C}\setminus \CC^{S} \longrightarrow \mathbb{C}^{2\times 2}
$ satisfies the following Riemann-Hilbert problem:
\begin{enumerate}
\item[$(RH6.1)$] $S$ is analytic in $\mathbb{C}\setminus \CC^{S}$.
\item[$(RH6.2)$] $S$ has continuous boundary values at regular
points of $\CC^{S}$, denoted by $S_+$ and $S_-$, such that $S_+(z)
= S_-(z)J_{S}(z)$, where
\[ J_{S}(z)= \begin{cases}
    P^{(j)}(z)N(z)^{-1},\, & z \in \partial
    \Delta_{\varepsilon}(\zeta_j) \,, \quad j=1, 2\,; \\
    N(z)J_{T}(z)N(z)^{-1},\, & z \in \CC^{S}\setminus (\partial \Delta_{\varepsilon}(\zeta_1)\cup \partial
    \Delta_{\varepsilon}(\zeta_2)).
    \end{cases}\]
\item[$(RH6.3)$] $S(z)\,=\,I+O\left(1/z\right),\,\,z\rightarrow
\infty$.
  \item[$(RH6.4)$] $S(z)$ is bounded in a neighborhood of the singular points of $\CC^S$.
\end{enumerate}
Observe that by construction, $J_S=I+O(1/n)$ as $n\to \infty$ on
$\partial \Delta_{\varepsilon}(\zeta_j)$, $j=1, 2$, and is
exponentially close to $I$ on the rest of contours of $\CC^S$.
Using the same arguments as in \cite{Deift99} we conclude that
\begin{equation*}\label{limit for V}
    S(z)= I+O\left(\frac{1}{n}\right)
\end{equation*}
uniformly for $z\in \mathbb{C}\setminus \CC^{S}$.

\section{Proofs of the main results}

We establish strong asymptotics for $\{p_n\}$ tracing back all the
previous transformations. For the sake of brevity, we do it
explicitly only for the case C.3. The proofs in the case C.2 are
very similar.

\subsection{Proof of Theorems \ref{theoremstrongOutside} and \ref{theoremstrongInside} in the case C.3}

Assume  that $-1 < A< 0 < B $ (case C.3); by \eqref{first transf},
$$
Y(z)=d_n ^{\sigma_3/2}\kappa^{-n\sigma_3}U(z) G(z)^{n\sigma_3} d_n
^{-\sigma_3/2},
$$
so that
\begin{equation}\label{Y11bis}
    p_{n}(z)= Y_{11}(z)= \left(\frac{G(z)}{\kappa}\right)^n\, U_{11}(z).
\end{equation}
Assume that $z\in \C \setminus \Pc(\Sigma)$, away from the branch
points; without loss of generality we may take $z\in
\Omega_\infty^T\setminus (\overline{\Delta_{\varepsilon}(\zeta_1)}
\cup   \overline{\Delta_{\varepsilon}(\zeta_2)})$ (see Figure
\ref{fig:FirstLensesCases23} or \ref{fig:FinalLensesCases23}).
Then by \eqref{T},
$$
U(z)= S(z)N(z)\,,
$$
and taking into account the expression of $N$ in (\ref{model}), we
obtain that uniformly on compact subsets of $\Omega_{\infty}^{T}
\setminus (\overline{\Delta_{\varepsilon}(\zeta_1)} \cup
\overline{\Delta_{\varepsilon}(\zeta_2)})$,
\begin{equation}\label{Y11}
    p_{n}(z)=Y_{11}(z)= \left(\frac{G(z)}{\kappa}\right)^n\,
    \left(S N \right)_{11}(z)=
    \left(\frac{G(z)}{\kappa}\right)^n\,N_{11}(z)\,\left(1+O\left(\frac{1}{n}\right)\right)\,,
\end{equation}
which proves \eqref{AsympOutside}.

If  $z$ belongs to the bounded component of $\C\setminus \Gamma$,
we may assume without loss of generality that $z\in
\Omega_1^T\setminus \overline{\Delta_{\varepsilon} (\zeta_2)}$. By
\eqref{T} and \eqref{V},
$$
U(z)=S(z) N(z) \begin{pmatrix} H^{-2n}(z) & -d_n^{-1} \\ d_n &
  0   \end{pmatrix} \,,
$$
and again uniformly in compact subsets of $\Omega_{1}^T \setminus
\Delta_{\varepsilon} (\zeta_2)$,
\begin{equation*}
\begin{split}
p_{n}(z)=& \left(\frac{G(z)}{ \kappa}\right)^n\, \left( \left[ S N
\right]_{11}(z) H^{-2n}(z)  + d_n \left[ S N \right]_{12}(z)
\right)
\\
=& \left(\frac{G(z)}{ \kappa}\right)^n\, \left(H^{-2n}(z)
 N_{11}(z) \left(1+O\left(\frac{1}{n}\right)\right)
+\,d_{n}N_{12}(z)
\left(1+O\left(\frac{1}{n}\right)\right)\right)\,,
\end{split}
\end{equation*}
which proves \eqref{AsympC3Gamma}. Obviously, this formula is
valid also if $z\in \Gamma_-$, that is, if $z$ lies on the
``$-$''-side of $\Gamma$, away from the branch points.

Assume now that $z\in \Gamma_+$ away from  $\zeta _2$. Again,
without loss of generality we may take $z\in \Omega_R \setminus
\overline{\Delta_{\varepsilon} (\zeta_2)}$. By \eqref{T} and
\eqref{V},
$$
U(z)=S(z) N(z) \begin{pmatrix} 1 & 0 \\
    d_n   H^{-2n}(z) & 1 \end{pmatrix}  \,,
$$
and uniformly in compact subsets of $\Omega_{R}  \setminus
\overline{\Delta_{\varepsilon}(\zeta_2)}$,
\begin{equation*}
\begin{split}
p_{n}(z)=& \left(\frac{G(z)}{ \kappa}\right)^n\, \left( \left[ S N
\right]_{11}(z)   + d_n H^{-2n}(z) \left[ S N \right]_{12}(z)
\right)
\\
=& \left(\frac{G(z)}{ \kappa}\right)^n\, \left(
 N_{11}(z) \left(1+O\left(\frac{1}{n}\right)\right)
+\,d_{n}H^{-2n}(z)
N_{12}(z)\left(1+O\left(\frac{1}{n}\right)\right) \right)\,,
\end{split}
\end{equation*}
which proves \eqref{AsympC3GammaBis}. Using
\eqref{boundaryValuesH} it is easy to see that formulas
\eqref{AsympC3Gamma} and \eqref{AsympC3GammaBis} match on
$\Gamma$.

Finally, if $z$ lies on the $\pm$-side of the interval $(\zeta _1,
\zeta _2)$, we assume $z \in \Omega_L^\pm\setminus
(\overline{\Delta_{\varepsilon}(\zeta_1)} \cup
\overline{\Delta_{\varepsilon}(\zeta_2)})$, where by \eqref{T} and
\eqref{V},
$$
U(z)=S(z) N(z) \cdot \begin{cases} \begin{pmatrix} 1 & 0 \\
       H^{-2n}(z) & 1 \end{pmatrix}  \,, & z\in \Omega_L^+\,, \\
     \begin{pmatrix} 1 & 0 \\
     -e^{-2A\pi i n} H^{-2n}(z) & 1 \end{pmatrix}  \,, & z\in \Omega_L^-\,.
\end{cases}
$$
Hence, uniformly in compact subsets of $\Omega_L^+\setminus
(\overline{\Delta_{\varepsilon}(\zeta_1)} \cup
\overline{\Delta_{\varepsilon}(\zeta_2)})$,
\begin{equation*}
\begin{split}
p_{n}(z)=& \left(\frac{G(z)}{ \kappa}\right)^n\, \left( \left[ S N
\right]_{11}(z)   +  H^{-2n}(z) \left[ S N \right]_{12}(z) \right)
\\
=& \left(\frac{G(z)}{ \kappa}\right)^n\, \left(
 N_{11}(z) \left(1+O\left(\frac{1}{n}\right)\right)
+H^{-2n}(z) N_{12}(z)\left(1+O\left(\frac{1}{n}\right)\right)
\right)\,,
\end{split}
\end{equation*}
while uniformly in compact subsets of $\Omega_L^-\setminus
(\overline{\Delta_{\varepsilon}(\zeta_1)} \cup
\overline{\Delta_{\varepsilon}(\zeta_2)})$,
\begin{equation*}
\begin{split}
p_{n}(z)=& \left(\frac{G(z)}{ \kappa}\right)^n\, \left( \left[ S N
\right]_{11}(z)   -e^{-2A\pi i n}  H^{-2n}(z) \left[ S N
\right]_{12}(z) \right)
\\
=& \left(\frac{G(z)}{ \kappa}\right)^n\, \left(
 N_{11}(z) \left(1+O\left(\frac{1}{n}\right)\right)
-e^{-2A\pi i n} H^{-2n}(z)
N_{12}(z)\left(1+O\left(\frac{1}{n}\right)\right) \right)\,.
\end{split}
\end{equation*}
This finishes the proof of
(\ref{AsympC3int})--\eqref{AsympC3intDown}.

\subsection{Proof of Theorem \ref{theoremstrongC3endpoint}}

Let $z \in \Delta_{\varepsilon} (\zeta_2)$; then by \eqref{V},
$$
U(z)=S(z) P(z) K^{-1}(z)\,,
$$
where $K(z)$ is one of the matrices given in the right hand side
of \eqref{T}. Gathering \eqref{Pcommon}--\eqref{En},  we get that
$$
P(z)=\sqrt{\pi}e^{\frac{\pi i}{6}}\,
    \begin{pmatrix}
    1 & -1 \\
    -i & -i \\
    \end{pmatrix}\,
    \left( \frac{t_n^{1/4} }{a(z)} \right)^{\sigma_3}\, \Psi \left( t_n \right)\,
     e^{-A \pi i n \sigma_3}\, e^{n\phi(z)
\sigma_3}\,,
$$
with $t_n\isdef n^{2/3} f(z)$ and $\Psi$ given by \eqref{Airy}.
For instance, if $z \in \Delta_{\varepsilon} (\zeta_2) \cap
\Omega_{R}$,  using \eqref{relationPhiH} we get
\begin{equation*}
\begin{split}
U(z)=& \sqrt{\pi}e^{\frac{\pi i}{6}}\,S(z)\,
    \begin{pmatrix}
    1 & -1 \\
    -i & -i \\
    \end{pmatrix}\,
    \left( \frac{t_n^{1/4}}{a(z)} \right)^{\sigma_3}\, \Psi \left( n^{2/3} f(z)\right)\\
     &\times
     e^{-A \pi i n \sigma_3}\, e^{n\phi(z)
\sigma_3}\,\begin{pmatrix}
    1 & 0 \\
    d_n e^{2 n\phi(z)} & 1 \\
    \end{pmatrix}\,.
\end{split}
\end{equation*}
Observe that for $z \in \Delta_{\varepsilon} (\zeta_2) \cap
\Omega_{R}^+$, we have $ 0<\arg f(z)<2\pi /3$, and we use the
expression
$$
\Psi(t)=
\begin{pmatrix}
    \Ai(t) & \Ai(\omega^{2}t) \\
    \Ai'(t) & \omega^{2}\Ai'(\omega^{2}t) \\
    \end{pmatrix} e^{-\frac{\pi i}{6} \sigma_{3}}\,.
$$
Hence,
\begin{equation*}
\begin{split}
\begin{pmatrix}
U_{11}(z) \\ U_{21}(z)
\end{pmatrix}=& \sqrt{\pi}\, e^{  n\phi(z)}\,S(z)\,
    \begin{pmatrix}
    1 & -1 \\
    -i & -i \\
    \end{pmatrix}\,
    \left( \frac{t_n^{1/4} }{a(z)} \right)^{\sigma_3}\, \begin{pmatrix}
   \mathcal A (t_n)   \\
   \mathcal A' (t_n)  \\
    \end{pmatrix}\,,
\end{split}
\end{equation*}
where
\begin{equation*}\label{defA}
    \mathcal A(t)\isdef e^{-A\pi i n} \, \Ai(t) +
2i \, e^{\frac{\pi i}{3}}\, \sin(A\pi n) \, \Ai (\omega^2
t)\,,\quad \omega=e^{2\pi i /3}\,.
\end{equation*}
Consequently,
\begin{equation*}\label{U11}
U_{11}(z)=\sqrt{\pi}\, e^{  n\phi(z)} \left(
\frac{t_n^{1/4}}{a(z)}\, \mathcal A (t_n) \, \left( 1+O\left(
\frac{1}{n}\right) \right)- \frac{a(z)}{t_n^{1/4}}\, \mathcal A'
(t_n) \, \left( 1+O\left( \frac{1}{n}\right) \right) \right)\,,
\end{equation*}
and taking into account \eqref{Y11bis} and the fact that in $
\Delta_{\varepsilon} (\zeta_2) \cap \Omega_{R}^+$,
$\exp(-\phi)=H$, we arrive at \eqref{localCaseC3}. Proceeding in a
similar way, we see that this expression is  also valid for $z$ in
the other regions of $ \Delta_{\varepsilon} (\zeta_2) $.

\subsection{Proof of Theorem \ref{theoremweak}}

This theorem is a corollary of Theorems \ref{theoremstrongOutside}
and  \ref{theoremstrongInside}. First, taking into account
\eqref{AsympOutside} and that function $N_{11}$, defined in
\eqref{defNN}, has no zeros in the plane cut from $\zeta _1$ to
$\zeta _2$, we see that zeros of $\{p_n\}$ cannot accumulate at
$\C\setminus \Pc (\Sigma)$.

Consider in particular case C.3. Now the asymptotic location of
the zeros of Jacobi polynomials depends  also on the value
\begin{equation*}\label{L}
   e^{-r }= \lim_{n\rightarrow \infty}|d_n|^{1/n}=
    \lim_{n\rightarrow \infty}|\sin (A\pi n)|^{1/n},
\end{equation*}
(assuming it exists), where $d_{n}$, defined in \eqref{defDn},
depends upon the distance of $\alpha_{n}=An$ to the integers, in
such a way that
$$
e^{-r } = \lim_{n\rightarrow \infty}(\dist
(\alpha_{n},\mathbb{Z}))^{1/n}.
$$
Let $z\in \Int(\Pc(\Sigma))$, that is, $z$ lies in the bounded
component limited by the contour $\Gamma$. We can choose
$\varepsilon  >0$ small enough such that $z \notin
\Delta_{\varepsilon}(\zeta_2)$. Then, the asymptotic formula
\eqref{AsympC3Gamma},
\begin{equation*}
\begin{split}
 p_n
(z)= &
\frac{1}{\kappa^n}\,\left(\left(G(z)w^2(z)\right)^{-n}N_{11}(z)
        \left(1 + O\left(\frac{1}{n}\right)\right)\right. \\ &
        \left. +2ie^{-An \pi i} \sin (A\pi n)\, G^n(z) N_{12}(z)
        \left(1 + O\left(\frac{1}{n}\right)\right) \right)\,,
\end{split}
\end{equation*}
is valid. Hence, $z$ is a zero of $p_n$ only if
$$
H^{-2n}(z)  =- 2ie^{-An \pi i} \sin (A\pi n)\, \frac{
N_{12}(z)}{N_{11}(z)}\, \left(1 +
O\left(\frac{1}{n}\right)\right)\,.
$$
Since $N_{12}/N_{11}$ is uniformly bounded and uniformly bounded
away from zero, we see that the zeros in this domain must satisfy
$$
|H^{-1}(z)|=| \sin (A\pi n)|^{1/(2n)}\, \left(1 +
O\left(\frac{1}{n}\right)\right)=e^{-r/2 }\, \left(1 +
O\left(\frac{1}{n}\right)\right)\,.
$$
It remains to use that $|H(z)|=e^{r/2 }$ defines in $
\Int(\Pc(\Sigma))$ the curve $\Gamma_r$.

Once we have established where the zeros accumulate, it remains to
prove that they asymptotically distribute according to the
corresponding measures in parts (i)-(ii). To this end we can use
the second order linear differential equation satisfied by Jacobi
polynomials $y_n =P_n ^{(\alpha_n,\beta_n)}$ (see e.g.\ \cite[\S
4.22]{szego:1975}):
$$
(1-z^2)\, y_n''(z) +\left[ \beta _n-\alpha _n-(\alpha _n+\beta
_n+2) z \right]\, y_n'(z) + n (n+\alpha _n+\beta _n+1)\,
y_n(z)=0\,.
$$
If we rewrite it in terms of $h_n=y'_n/(n y_n)$, we obtain a
Riccati differential equation:
\begin{equation}\label{riccati}
(1-z^2)\, \left( \frac{1}{n}\, h_n'(z)+h_n^2(z) \right) +
\frac{\beta _n-\alpha _n-(\alpha _n+\beta _n+2) z }{n} \, h_n'(z)
+ \frac{ n+\alpha _n+\beta _n+1 }{n}=0\,.
\end{equation}
Let $\nu_n$ denote the normalized zero counting measures of
$y_n=P_n^{(\alpha _n, \beta _n)}$. By a weak compactness argument
we know that there exists an infinite subsequence $\Lambda \subset
\mathbb{N}$ and a unit measure $\nu$ such that $\nu_n \rightarrow
\nu,\,n\in \Lambda$, in the weak*-topology. In the first part of
this proof, we saw that $\supp(\nu) $ consists of a finite union
of analytic arcs or curves, and every compact subset of
$\C\setminus \supp (\nu)$ contains no zeros of $P_n^{(\alpha _n,
\beta _n)}$ for $n$ sufficiently large.

Hence,
$$
h_n(z)=\int_\Gamma \frac{d\nu_n(t)}{z-t} \longrightarrow
h(z)=\int_\Gamma \frac{d\nu(t)}{z-t}\,, \quad n\in \Lambda\,,
$$
locally uniformly in $\C\setminus \supp (\nu)$. Taking limits in
(\ref{riccati}) we obtain that $h$ satisfies the quadratic
equation
$$
(1-z^2)\, h^2(z)  + \left[ B-A-(A+B) z \right]\, h(z) + A +B+1
=0\,,
$$
so that
$$
\int_\Gamma \frac{d\nu(t)}{z-t}=\frac{A+B+2}{2}\,
\frac{R(z)}{z^2-1}-\frac{1}{2}\, \left( \frac{A}{z-1}+
\frac{B}{z+1}\right)\,, \quad z \in \C\setminus \supp (\nu)\,.
$$
By Sokhotsky-Plemelj's formulas, on every arc of $\supp (\nu)$,
\begin{equation}\label{nu}
   d\nu(z) = \frac{A + B+2}{2 \pi i} \frac{R_+(z)}{z^2-1} \,dz\,.
\end{equation}

Now, we are concerned with proving part (ii) of the theorem,
related to the case C.3. First, consider the generic case when
$r=0$. In this case, the measure $\mu$ in (\ref{measure}) is
supported on $\Sigma =\Gamma \cup [\zeta_1,\zeta_2]$. Thus, by
\eqref{nu}, $\mu'=\nu'$ a.e.\ on $\supp(\nu)$, $\mu ,\nu$ being
unit measures. Therefore, $\nu =\mu$. The proof in the case
$0<r<\infty$ is similar, but with measures $\mu_r$, given in
(\ref{measure level}), in place of $\mu$. Finally, for the
degenerate case $r=\infty$, which takes place when parameters
$\alpha_n$ approach the integers faster than exponentially, it is
enough to take into account that the Cauchy transform of the
measure $d\sigma =-A\delta_1$ is $\widehat{\sigma} (z)=-A/(z-1)$.

Finally, the proof for the Case C.2 is totally analogous.

\section*{Acknowledgements}

The research of A.M.F.\ was supported, in part, by a research
grant from the Ministry of Science and Technology (MCYT) of Spain,
project code BFM2001-3878-C02, by NATO Collaborative Linkage Grant
``Orthogonal Polynomials: Theory, Applications and
Generalizations'', ref. PST.CLG.979738, and by Research Network
``Network on Constructive Complex Approximation (NeCCA)'', INTAS
03-51-6637.

The research of R.O.\ was partially supported by grants from
Spanish MCYT (BFM2001-3411) and  Gobierno Aut\'{o}nomo de Canarias
(PI2002/136).

Finally, the authors are grateful to Prof. A. B. J. Kuijlaars for
stimulating discussions on the Riemann-Hilbert analysis.

\end{document}